\documentclass[twocolumn]{autart}
\usepackage{cite}
\usepackage{graphicx}
\usepackage{color}
\usepackage{textcomp}
\usepackage{amsmath}
\usepackage{booktabs}
\usepackage{amsfonts}
\usepackage{amssymb}
\usepackage{mathtools}
\usepackage{bm}
\usepackage{float}
\usepackage{comment}
\usepackage{subcaption}
\usepackage[normalem]{ulem}
\usepackage{multirow}
\usepackage[font={small}]{caption}
\usepackage{tabstackengine} 
\usepackage[mathscr]{euscript}
\usepackage{ulem}
\usepackage{bbm}
\graphicspath{ {images/} }
\usepackage[dvipsnames]{xcolor}

\usepackage{enumitem}

\begin{document}
\begin{frontmatter}
 \title{The Effect of Sensor Fusion on Data-Driven Learning of Koopman Operators}
\author[UCSB-EE]{Shara Balakrishnan}\ead{sbalakrishnan@ucsb.edu},    
\author[UCSB-ME]{Aqib Hasnain}\ead{aqib@ucsb.edu},               
\author[PNNL]{Rob Egbert}\ead{rob.egbert@pnnl.gov},  
\author[UCSB-ME]{Enoch Yeung}\ead{eyeung@ucsb.edu}              
\address[UCSB-EE]{Department of Electrical and Computer Engineering, University of California, Santa Barbara, United States} 
\address[UCSB-ME]{Department of Mechanical Engineering, University of California, Santa Barbara, United States} 
\address[PNNL]{Biological Sciences Division, Earth and Biological Sciences Directorate, Pacific Northwest National Laboratory, United States} 
\begin{abstract} 
Dictionary methods for system identification typically rely on one set of measurements to learn governing dynamics of a system. In this paper, we investigate how fusion of output measurements with state measurements affects the dictionary selection process in Koopman operator learning problems. While prior methods use dynamical conjugacy to show a direct link between Koopman eigenfunctions in two distinct data spaces (measurement channels), we explore the specific case where output measurements are nonlinear, non-invertible functions of the system state. This setup reflects the measurement constraints of many classes of physical systems, e.g., biological measurement data, where one type of measurement does not directly transform to another. We propose output constrained Koopman operators (OC-KOs) as a new framework to fuse two measurement sets. We show that OC-KOs are effective for sensor fusion by proving that when learning a Koopman operator, output measurement functions serve to constrain the space of potential Koopman observables and their eigenfunctions. Further, low-dimensional output measurements can be embedded to inform selection of Koopman dictionary functions for high-dimensional models.  We propose two algorithms to identify OC-KO representations directly from data: a direct optimization method that uses state and output data simultaneously and a sequential optimization method. We prove a theorem to show that the solution spaces of the two optimization problems are equivalent. We illustrate these findings with a theoretical example and two numerical simulations.

\end{abstract}
\end{frontmatter}


\section{Introduction}



The rapid evolution in sensor technology enables us to acquire a wealth of information about the governing dynamics of nonlinear systems. Novel sensors are being designed and fabricated at an rapid rate to enhance the data acquisition process for systems in various domains like biology \cite{kim2019wearable,gupta2019cell,bhalla2020opportunities,justino2017graphene,metkar2019diagnostic}, mechanics \cite{bhatt2019mems,abels2017nitride,tao2012gait, fiorillo2018theory, taheri2019review}, and transportation \cite{guerrero2018sensor, bachmann2011multi, richoz2020transportation, zhao2019detection, bernas2018survey}.
The tools, techniques, and theories that integrate all the data to augment our knowledge of the system is the broad area of sensor fusion \cite{mitchell2007multi}.
Sensor fusion plays a pivotal role in many applications.
In the human gait system, senor fusion methods are deployed to study gait dynamics \cite{meng2010emg,ziegler2018classification,zhao2019adaptive, horst2019explaining,alharthi2019deep}, detect gait anomalies \cite{ajay2018pervasive, paragliola2018gait, bonetto2019seq2seq} and control prosthetics \cite{grimmer2016powered,guo2010study,hoover2012stair,au2005emg,huang2011continuous,zhang2019sensor}. Numerous sensors are utilized in inertial navigation systems with data fusion architectures for accurate state estimation to aid in navigation and control \cite{sasiadek2002sensor,hasan2009review,hansen2017nonlinear,shankar2019finite,rigatos2007extended,rigatos2010extended,yazdkhasti2018multi,hegrenaes2007towards}. Maximizing the throughput of manufacturing processes \cite{shi2018using,pratama2019metacognitive,dou2020unsupervised,kuntouglu2021investigation,bleakie2013feature,lin2021metal,francis2019deep,shevchik2017acoustic,gu2017product,he2017integrated}, minimizing traffic congestion \cite{bachmann2011multi,richoz2020transportation,zhao2019detection}, and monitoring structural health \cite{wu2020data,vitola2017sensor,ostachowicz2019optimization} are some of the other areas where sensor fusion plays a critical role.


In recent years, there has been significant interest in developing sensor fusion techniques for biological systems, spurred in part by decreasing costs of next-generation omics measurements
\cite{huang2017more,dai2021multi, rohart2017mixomics,noor2019biological}. While some techniques like transcriptomics \cite{wang2009rna} and proteomics \cite{aebersold2003mass} inform the genetic activity within the cell, other instruments like flow cytometers \cite{gant1993application}, plate readers \cite{meyers2018direct}, and microscopes \cite{lefman2004three} inform the phenotypic characteristics outside the cell. Substantial progress towards fusing a combination of these datasets with prior knowledge of the system has resulted in static models that provide insights about the underlying network topology of the interaction between the genes and the proteins \cite{huang2017more, dai2021multi, altemose2020mudamid, liu2020integrating, rohart2017mixomics,noor2019biological, hwang2005data_simulation, hwang2005data}. 

An area of research emerging in the last decade is the problem of sensor fusion for biological networks, that couple distinct streams of omics measurements with fluorescence data.  Fluorescence data is often used to parameterize or identify dynamical system models describing governing dynamics and network topology, while omics measurements are used to deduce whole-cell statistics and steady-state phenotypes of cellular metabolism, stress, and fitness.   To that end, our goal in this paper is to develop sensor fusion techniques that integrate various types of time-series data to construct dynamic genotype-to-phenotype models.

Koopman operator methods have recently demonstrated great promise in simultaneously discovering A) the governing dynamics and B) a spectral decomposition of a complex physical system represented purely by data \cite{mezic2005spectral,budivsic2012applied}. The key premise of Koopman operator theory is that a collection of state functions, or observables, can be constructed, discovered, or estimated, to embed nonlinear dynamics of a physical system  in a high dimensional space. The Koopman operator then acts as a linear operator on the function space of observables, governing time-evolution of the system dynamics, as a linear system.    

In reality, the Koopman operator is infinite dimensional and must be approximated numerically.   The problem of finding a Koopman operator and a collection of observable functions is known as the Koopman operator learning problem. 
The classical approach to solve this problem is to use dynamic mode decomposition (DMD) \cite{schmid2010dynamic}.  More recently, variations on DMD involve approximating observable functions with a broad set of dictionary functions (E-DMD) \cite{williams2015data,hasnain2019data}, which can be generated using deep learning \cite{lusch2018deep,yeung2019learning,kaneko2019convolutional,azencot2019consistent}, by casting the learning problem, as a robust optimization problem to handle sparse data \cite{askham2018variable,sinha2019computation} or to treat heterogeneously sampled data \cite{manohar2018optimal,manohar2019optimized}. 
The power of Koopman operators lies in their ability to capture the underlying modes that drive the system \cite{mezic2013analysis,taira2017modal,mclean2020modal}, directly from data. Koopman operators also enable the construction of observers \cite{surana2016linear,surana2016koopman,yeung2018koopmanGramian,netto2018robust} and controllers \cite{proctor2016dynamic,korda2018linear,proctor2018generalizing,you2018deep,kaiser2021data} for nonlinear systems in a linear framework.

Koopman operators can fuse sensor measurements into a single dataset to learn the underlying dynamics, as seen in systems like traffic dynamics \cite{avila2020data,ling2020koopman}, human gait \cite{boudali2017human,kalinowska2019data, fujii2019data}, and robotics \cite{haggerty2020modeling}.
The work of Williams et al. \cite{williams2015datafusion} and Mezic \cite{mezic2019spectrum} use Koopman operator theory to elucidate the fusion of the dynamics evolving on two different state-spaces, provided that there is a function map between the state-spaces. Williams et al. \cite{williams2015datafusion} consider two datasets that are rich enough to reconstruct the system state and develop an algorithm to map the eigenfunctions of the Koopman operator identified from each dataset.
Mezic \cite{mezic2019spectrum} proves a relational mapping between eigenfunctions of both spaces, when exact conjugacy is not possible. For dynamics that evolve on two different state-spaces, they define the factor conjugacy of the dynamics for the function that maps the two spaces. 

 This paper builds on existing Koopman operator fusion theory \cite{williams2015datafusion,mezic2019spectrum} to examine a special case: we consider learning a sensor fusion model for a physical system represented by direct state data and a series of output measurements.  We consider the scenario where both the Koopman operator, observables, and the relational map between Koopman observable and output measurements are unknown.  In a standard Koopman operator learning problem, the state measurement data is sufficient to approximate the Koopman operator.  However, we examine the effect of {\it adding} output measurements now as a series of behavioral constraints on the dynamics of the Koopman operator --- we aim to know the effect of incorporating output measurements (sensor fusion) on the solution of the Koopman operator learning problem.  Specifically, we seek to understand how Koopman eigenfunctions, spectra, and modes change as a consequence of sensor fusion.

The formulation of an output-constrained Koopman operator is not novel.  In the literature, output-constrained Koopman operators are used for various applications like observability analysis \cite{mesbahi2021nonlinear,yeung2018koopmanGramian}, observer synthesis \cite{surana2016linear,surana2016koopman,netto2018robust}, and sensor placement \cite{manohar2019optimized,manohar2018optimal,hasnain2019optimalReporterPlacement} for nonlinear systems.  In this paper, we prove that output-constrained Koopman operators satisfy the following properties:
\begin{enumerate}[label=(\textit{\roman*})]
    \item The output dynamics of the nonlinear system always span a subspace of observable functions for the output-constrained Koopman operator \vspace{2pt}
    \item The observables of the output-constrained Koopman operators can capture the dynamics of both states and outputs \vspace{2pt}
    \item State-inclusive output-constrained Koopman operators exist in the region of convergence of the Taylor series expansion of the dynamics and output functions of any nonlinear system.
\end{enumerate}
Heretofore, there have been few algorithms that identify output-constrained Koopman operators, to the best of our knowledge. 
To identify output-constrained Koopman operators (OC-KO) from data, we pose the output-constrained DMD (OC-DMD) problem as a special extension of the DMD problem to incorporate output constraints. We propose two variants of the problem: the \textit{direct OC-DMD} solves for the state and output dynamics concurrently, while the \textit{sequential OC-DMD} solves for them sequentially. 
Sequential OC-DMD explicitly reveals the effect of having outputs in the KOR learning problem. To implement OC-DMD in practice, we build on the deepDMD algorithm developed by Yeung et al. \cite{yeung2019learning}, where neural networks represent vector valued observables of the Koopman operator. We then study the effect of affine transformations on the output-constrained Koopman operator learning problem to take into account--- how data preprocessing methods like normalization or standardization modify the output-constrained Koopman operator. We use simulation examples to investigate the performance of OC-DMD algorithms. Our findings include:
\begin{enumerate}[label=(\textit{\roman*})]
    \item The solution space of the direct OC-DMD and sequential OC-DMD optimization problems are equivalent \vspace{5pt}
    \item Affine state transformations yield OC-KOs with an eigenvalue on the unit circle \vspace{5pt}
    \item OC-KOs optimized for multi-step predictions are required to capture the dynamics with limit cycles.
\end{enumerate}

The paper is organized as follows. In section \ref{sec: Basic Idea}, we formulate the OC-KO representation. In section \ref{sec: math prelims}, we briefly discuss Koopman operator theory, Koopman operator sensor fusion  and the necessary DMD algorithms. We discuss the properties of OC-KOs in Section \ref{sec: oc- Koopman} and the OC-DMD algorithms in Section \ref{sec: oc-deepDMD}. We show the simulation results in Section \ref{sec: Simulation Results} and conclude our analysis in Section \ref{sec: Conclusion}.

\section{Problem Formulation} \label{sec: Basic Idea}
Our goal in this paper will be to consider physical systems represented by sampled time-series data, even their underlying governing dynamics are continuous.  For methods that estimate the Koopman generator (the continuous-time extension of the discrete-time Koopman operator), we refer the reader to \cite{mauroy2020koopman}.  Suppose we have an autonomous discrete-time nonlinear dynamical system with output
\begin{subequations}\label{eq: nonlinear system with output}
    \begin{flalign}
        \text{State Equation:} && x_{k+1} &= f(x_k)&& \label{eq: NL_sys_without_output}\\
        \text{Output Equation:} && y_k &= h(x_k)&& \label{eq: NL_output_equation}
    \end{flalign}
\end{subequations}
where $x \in \mathcal{M} \subseteq \mathcal{R}^n$ is the state, $y \in \mathbb{R}^p$ is the output, $f:\mathcal{M} \rightarrow \mathcal{M}$ and $h: \mathcal{M} \rightarrow \mathbb{R}^p$ are analytic functions and $k$ is the discrete time index indicating the time point $kT_s$ with $T_s$ being the sampling time.  Let 
\begin{equation}\label{eq: Koopman system with output}
    \begin{aligned}
        \psi(x_{k+1}) &= K\psi(x_{k})\\
        y_{k} &= W_h\psi(x_{k})
    \end{aligned}
\end{equation}
be the OC-KO  of (\ref{eq: nonlinear system with output}) where $\psi:\mathcal{M} \rightarrow \mathbb{R}^{n_L}$ is a vector function of state-dependent scalar observable functions ( $n_L \leq \infty$), $K \in \mathbb{R}^{n_L\times n_L} $ is the Koopman operator and $W_h \in \mathbb{R}^{p\times n_L} $ is a projection matrix that projects observables to the space of output functions. We record state and output measurements ($x$ and $y$, respectively) from (\ref{eq: nonlinear system with output}) as
\begin{equation}\label{eq: data matrices x}
    \begin{aligned}
        X_P &= [x_0,\cdots,x_{N-1}], \quad X_F = [x_1, \cdots,x_N],\\
        Y_P &= [y_0,\cdots,y_{N-1}],
    \end{aligned}
\end{equation}
where $X_P \in \mathbb{R}^{n\times N}$ is the state data collected from time points $0,...,N-1$, and  $X_F \in \mathbb{R}^{n\times N}$  and $Y_P \in \mathbb{R}^{p\times N}$ are the collection of state and output measurements propagated one step forward from the elements in $X_P$, respectively. Here, we use lower-case notation for the state $x$ and output $y$ variables and upper-case notation for variables containing sampled time-series snapshots involving the state variable $x$ with $X_P,X_F$ and time-series with output variable $y$ with $Y_P$, respectively. 

\section{Mathematical Preliminaries}\label{sec: math prelims}
We now briefly introduce the formal mathematical elements of Koopman operator theory, its modal decomposition, Koopman operator sensor fusion and relevant DMD algorithms. 

\subsection{Koopman operator}
Given the dynamical system (\ref{eq: NL_sys_without_output}), the KO represented by $\mathcal{K}$ is a linear operator that is invariant in the functional space $\mathcal{F}$ as $\mathcal{K}: \mathcal{F} \rightarrow \mathcal{F}$. Any function $\tilde{\phi} \in \mathcal{F}$ such that $\tilde{\phi}:\mathcal{M}\rightarrow \mathbb{C}$ is defined as a \textit{scalar observable} with the property
\[
(\mathcal{K}\tilde{\phi})(x) = (\tilde{\phi} \circ f)(x)
\]
where $\tilde{\phi} \circ f \in \mathcal{F}$ due to invariance of $\mathcal{K}$.  
Let the set of basis functions of the function space $\mathcal{F}$ be denoted by 
\begin{align} \label{eq: basis functions}
    \Phi \triangleq \{\phi_1,\phi_2,\hdots,\phi_M\} \text{ with } M \rightarrow \infty.
\end{align}
Then, any function $\tilde{\phi} \in \mathcal{F}$ can be written as a linear combination of the basis functions $\phi_i \in \Phi$ implying $\tilde{\phi} = a^T\Phi$ with $a \in \mathbb{R}^{1 \times M}$. A vector valued observable (also referred to as a dictionary of observables) can be constructed by taking a vector combination of the basis functions in $\mathcal{F}$:
\begin{align*}
    \varphi(x) &= a^T\Phi(x) = \sum_{i=0}^M a_i\phi_i(x)
\end{align*}
where $\varphi:\mathcal{M}\rightarrow \mathbb{C}^{n_{\varphi}}$, $a_i \in \mathbb{R}^M$, $b_i \in \mathbb{R}^{n_{\varphi}}$,

Then, $\varphi(x)$ is invariant under $\mathcal{K}$ as
\begin{align*}
    \mathcal{K}\varphi(x_k) & = (\mathcal{K}a^T\Phi)(x_k) = a^T(\mathcal{K}\Phi)(x_k)\\
    &=a^T(\Phi \circ f)(x_k)= (a^T\Phi)(x_{k+1}) = \varphi(x_{k+1}).
\end{align*}
Thus, every linear combination of basis functions in $\Phi$ is also a Koopman observable; specifically the Koopman operator is a linear operator on the function space $\mathcal{F}$ spanned by $\Phi$.  

These observations also hold, more generally, when ${\tilde{\phi}}$ is a vector-valued observable and $\mathcal{F}$ defines a space of vector-valued functions.   The Koopman operator $\mathcal{K}$, in this setting, acts on vector-valued functions as opposed to scalar-valued functions.

The choice of $\varphi$ should be such that the  true state $x$ can be recovered. One such choice of $\varphi$ is given by
\begin{equation} \label{eq: definition state inclusive observable function}
    \psi(x) =\begin{bmatrix}
    x^T & \varphi^T(x)
\end{bmatrix}^T
\end{equation}
where $\psi: \mathcal{M} \rightarrow \mathbb{R}^{n_L}$ contains the base states $x$ in addition to $\varphi$. Such an observable $\psi(x)$ is called a state-inclusive Koopman observable \cite{johnson2018class}.  For the rest of the paper, unless explicitly stated,  we denote a $\varphi$ to be a dictionary of nonlinear observables (not necessarily state-inclusive) and $\psi$ to be a dictionary of state-inclusive observables with the form (\ref{eq: definition state inclusive observable function}).

\subsection{Modal decomposition}
The Koopman operator is infinite-dimensional as it acts on the functional space $\mathcal{F}$. As presented in \cite{williams2015data}, we take the basis functions $\Phi$ in (\ref{eq: basis functions}) to also be the set of eigenfunctions for $\mathcal{K}$ with $\lambda_i$ being the eigenvalue of $\phi_i$ and
$
\mathcal{K}\phi_i= \lambda_i\phi_i \quad \forall i \in \{1,2,\hdots,M\} \quad M \rightarrow \infty.
$
In this work we will assume that every Koopman operator describes the dynamics of an analytical system represented as (\ref{eq: nonlinear system with output}).  Such systems admit Koopman operators with countable spectra.  In this setting, the modal decomposition of the Koopman operator dynamics becomes
\begin{align}\label{eq: Koopman operator dynamics only}
    \varphi(x_{k+1}) &= \varphi \circ f(x_k) = \mathcal{K}\varphi(x_k) \nonumber\\
    &= \sum_{i=1}^M b_i \mathcal{K}\phi_i(x_k) = \sum_{i=1}^M b_i \lambda_i\phi_i(x_k)
\end{align}
where $b_i \in \mathbb{R}^{n_\varphi}$ are called the Koopman modes, $\lambda_i$ are the Koopman eigenvalues and $\phi_i$ are the corresponding Koopman eigenfunctions.  

\subsection{Koopman operators for conjugate dynamical systems} \label{sec: Fusion using Koopman operators}
We review the theory of eigenfunction conjugacy for conjugate dynamical systems \cite{mezic2019spectrum}. This theoretical construction informs how we view sensor fusion of dynamical systems, at large. Consider two nonlinear systems
\[
z_{k+1}^{(i)} = f_i(z_{k+1}^{(i)}), \hspace{5pt}z^{(i)}\in \mathbb{R}^{n_i},\hspace{5pt}i=1,2
\]
with their corresponding Koopman operators  $\mathcal{K}_{i}$. The two dynamical system are said to be factor conjugate if there is a function map $H:\mathbb{R}^{n_1} \rightarrow \mathbb{R}^{n_2}, n_1 \geq n_2$ where
\begin{align*}
    z^{(2)} &= H(z^{(1)})
    \text{ such that } H \circ f_1 (z^{(1)}) = f_2 \circ H (z^{(1)}).
\end{align*}
Then, if $\phi_{2}$ is an eigenfunction of the Koopman operator $\mathcal{K}_2$ corresponding to the eigenvalue $\lambda_{2}$, $
\mathcal{K}_2\phi_2 (z^{(2)}) = \lambda_2\phi_2 (z^{(2)})= \phi_2 \circ f_2 (z^{(2)})$. Then, we can see that $(\phi_2 \circ H)$ is an eigenfunction of $\mathcal{K}_1$:
\begin{align}\label{eq: sensor fusion KO eigenfunction map}
    \lambda_2 (\phi_2 \circ H)(z^{(1)}) \hspace{-2pt}&= \hspace{-2pt} \phi_2 \circ f_2 \circ H(z^{(1)}) \hspace{-2pt}=\hspace{-2pt} (\phi_2 \circ H) \circ f_1 (z^{(1)})\nonumber\\
    &= \mathcal{K}_1 (\phi_2 \circ H)(z^{(1)})
\end{align}
The dynamics of the second system spanned by a set of eigenfunctions can be mapped to the eigenfunctions of the first system. 
If $H$ is a $C^k$ diffeomorphism, i.e., $H$ has an inverse and is $k-times$ differentiable, then we have a $C^k$ diffeomorphic conjugacy. In that case, we have that for all the eigenfunctions $\phi_i$ such that $ \mathcal{K}_i\phi_i(z^{(1)})= \lambda \phi_i(z^{(i)}),i=1,2$ have a bijective map:
\begin{equation}\label{eq: diffeomorphic conjugacy eigenfunction map}
    \phi_1(z^{(1)}) = \phi_2 \circ H(z^{(1)}), \hspace{3pt} \phi_2(z^{(2)}) = \phi_1 \circ H^{-1}(z^{(2)}).
\end{equation}
The bijective map of eigenfunctions allows complete transfer of information between the two systems. We use this concept to identify subspaces of the state and output dynamics that are $C^k$ diffeomorphic conjugate and deduce how much information can be fused.

\subsection{Dynamic Mode Decomposition}

There are a class of nonlinear systems that have exact finite invariant Koopman operators \cite{brunton2016koopman}. 
But in most cases, the KOs are infinite-dimensional and difficult to identify numerically. The DMD algorithm introduced in \cite{schmid2010dynamic} finds KOs that are exact solutions for linear systems and local approximators for nonlinear systems. DMD uses the data matrices $X_p$ and $X_f$ as defined in (\ref{eq: data matrices x}) to solve
$
    K = \min_{K} ||X_f - K X_p|| = X_f X_p^{\dag}
$
where \textsuperscript{$\dag$} denotes the Moore-Penrose pseudoinverse. 

The extended DMD (E-DMD) algorithm proposed in \cite{williams2015data} identifies KOs that capture the nonlinear dynamics with better accuracy. E-DMD fuses kernel methods in machine learning with DMD to identify a rich set of observables to solve the optimization problem
\begin{equation*}
    K = \min_{K} ||\psi(X_f) - K \psi(X_p)|| =\psi(X_f)\psi(X_p)^{\dag}.
\end{equation*}
E-DMD solves the nonlinear regression problem using linear least squares. \cite{korda2018convergence} shows the relevance of E-DMD to Koopman operators. E-DMD hinges on the user specifying the lifting functions, and more often than not, it leads to an explosion of the lifting functions \cite{johnson2018class,balakrishnan2020prediction}.   

Recent developments in DMD incorporate deep learning approaches to identify the observables using deep neural networks. \cite{yeung2019learning,otto2019linearly,takeishi2017learning,li2017extended}. They can approximate exponentially many distinct observable functions. We consider the deep DMD formulation from \cite{yeung2019learning}:
\begin{align} \label{eq: deep DMD formulation}
\min_{\psi,K}&||\psi(x_{k+1})-K\psi(x_k)||^2_F\\
\psi(x) = 
\begin{bmatrix}
    x\\
    \varphi(x)
\end{bmatrix} =& 
\begin{bmatrix}
    x\\
    g_n \circ \sigma \circ \cdots \circ \sigma \circ g_2 \circ \sigma \circ g_1(x)
\end{bmatrix}\nonumber
\end{align}
where $\psi(x)$ is represented by neural network representations with the $i^{th}$ hidden layer captured by weights $W_i$, biases $b_i$, linear function $g_i(x) = W_ix+b_i$ and activation function $\sigma$. Hence, the estimation of $\psi$ boils down to learning the parameter set $(W_1,b_1,W_2,b_2,\hdots,W_n,b_n)$, while specifying or optimizing the depth of the network $n$ as a hyper-parameter. By selecting appropriate activation functions for the nonlinear transformation $\sigma$, in a given layer, e.g., using sigmoidal \cite{cybenko1989approximation}, rectified linear unit (ReLU) activation functions  \cite{hanin2019universal}  radial basis functions (RBFs) \cite{lo1998multilayer}, $\psi(x)$ leverages universal function approximation properties of each of these function classes \cite{yeung2019learning}. 

To derive the approximate Koopman eigenfunctions for any of these numerical approximations of the Koopman operator, suppose $KV =V\Lambda $ is the eigendecomposition of K where $v_i$ (the $i^{th}$ column of $V$) is the right eigenvector corresponding to the eigenvalue $\lambda_i$ (the $i^{th}$ diagonal entry of the diagonal matrix $\Lambda$). Then, $\psi(x_{k+1}) = K\psi(x_k)$ has the modal decomposition
\begin{equation} \label{eq: modal decomposition practical}
    \psi(x_{k+1}) = K\psi(x) = V\Lambda V^{-1}\psi(x) = \sum_i v_i\lambda_i \phi_i(x)
\end{equation}
where $\phi_i(x)$  is the $i^{th}$ entry of the vector $V^{-1}\psi(x)$. Comparing with (\ref{eq: Koopman operator dynamics only}), $\phi_i(x)$ and $v_i$ are the eigenfunction and the Koopman mode corresponding to the eigenvalue $\lambda_i$.  Thus, a spectral decomposition on any finite approximation $K$ of the true Koopman operator $\mathcal{K}$ gives an approximation to a subset of the true eigenfunctions of $\mathcal{K}.$

We remark that the true Koopman observables and true Koopman eigenfunctions for a given system often span an infinite dimensional space.  In the case of analytic state update equations for the system (\ref{eq: nonlinear system with output}), the dimension of the Koopman operator is usually countably infinite.  Only in special cases is exactly finite.  This means that any finite-dimensional set of functions may not exactly span or recover the Koopman observable function space.  In the sequel, we will refer to a finite collection of Koopman observable functions
$\psi_1(x), ...., \psi_{n_L}(x)$ as a dictionary of Koopman observables.  However, these are, in practice, an approximation of a true spanning set for the Koopman observable function space. 

\section{Output Constrained Koopman Operator Representations}\label{sec: oc- Koopman}

We consider the fusion of state and output measurements from the nonlinear system (\ref{eq: nonlinear system with output}) using the Koopman operator sensor fusion method as delineated in Section \ref{sec: Fusion using Koopman operators}. To do so, we need to establish a factor conjugate map between the state dynamics and the output dynamics. 

We suppose the state dynamics is captured by the KO $\psi(x_{k+1}) = K\psi(x_k)$. To integrate the state dynamics with the outputs (that are nonlinear functions of the state), we consider the OC-KO (\ref{eq: Koopman system with output}).

The model structure of OC-KO is such that its dictionary of observables  $\psi(x)$ span the nonlinear output functions $h(x)$ by augmenting the KO with a linear output equation $y_k = h(x_k) = W_h\psi(x_k)$. To establish a conjugacy between the state and output dynamics, all we need is to do is construct the dynamics of the output. The output dynamics can be given by
\begin{equation}\label{eq: output dynamic invertibility}
    y_{k+1} = W_h\psi(x_{k+1}) = W_hK\psi(x_k)= W_hKW_\psi y_k
\end{equation}
where $\psi(x_k)=W_\psi y_k$. We begin by considering the simplest case of inverting the $W_h$ matrix in $y_k = W_h\psi(x_k)$ to identify the $W_\psi$ matrix.

\begin{thm}\label{theorem: outputs span all lifting functions}
 Given a nonlinear system (\ref{eq: nonlinear system with output}) with a Koopman operator (\ref{eq: Koopman system with output}). Suppose $\psi(x) \in \mathbb{R}^{n_L}$ is the vector of Koopman observable functions for (\ref{eq: Koopman system with output}), not necessarily state-inclusive. 
For the number of outputs $p\geq n_L$, the following statements are true. 
\begin{enumerate}[label=(\textit{\roman*})]
    \item If $\quad rank(W_h)=n_L$, then each scalar Koopman observable function $\psi_i(x)$ lies in the span of the output functions $h(x)$ 
    \item If $\quad rank(W_h)=r<n_L$, then there exists a similarity transform $T$ that takes the model (\ref{eq: Koopman system with output}) to the form 
    \begin{equation}\label{eq: theorem 1 similarity transform}
        \begin{aligned}
            \Tilde{\psi}(x_{k+1}) &= \Tilde{K}\Tilde{\psi}(x_k)\\ 
            y_k &= \Tilde{W}_h\Tilde{\psi}(x_k)
        \end{aligned}
    \end{equation}
    where $\Tilde{K} = T^TKT$, $\Tilde{\psi}(x) = T^T\psi(x)$, $\Tilde{W}_h = W_hT$ such that $r$ components of the vector Koopman observable $\tilde{\psi}(x)$ given by $\psi_h(x)$ lie in the span of $h(x)$.
\end{enumerate}
\end{thm}

\begin{pf}
\textit{Case (1):} Since $W_h$ has full column rank, an exact left inverse exists such that 
\begin{align*}
    \psi(x_k) = (W_h^TW_h)^{-1}W_h^Th(x_k) \quad \forall \quad x_k \in \mathcal{M}
\end{align*}
and hence
$\psi_i \in span \{ h_1 , h_2, ..., h_p \}$ $\forall$ $i \in \{1,2,...,n_L\}$.\\
\textit{Case (2):} Suppose $rank(W_h) = r <n_L$. The singular value decomposition of $W_h$ yields
\begin{align} \label{eq: theorem 1 case 2}
W_h &= \begin{bmatrix}
        U_r & U_{p-r}
        \end{bmatrix}_{p\times p}
        \begin{bmatrix}
        \Sigma_r & 0\\
        0 & 0
        \end{bmatrix}_{p\times n_L}
        \begin{bmatrix}
        V_r^T\\ 
        V_{n_L-r}^T
        \end{bmatrix}_{n_L \times n_L}\nonumber\\
\Rightarrow
    h(x_k) &= \begin{bmatrix}
        W_{h\psi} & 0
        \end{bmatrix}
        \begin{bmatrix}
        \psi_h(x_k)\\ 
        \bar{\psi}_{h}(x_k)
        \end{bmatrix} \quad \forall \quad x_k \in \mathcal{M}\\
        \text{s. t. } W_{h\psi} \hspace{-3pt}&=\hspace{-3pt} U_r\Sigma_r,\hspace{1.5pt} \psi_h(x) \hspace{-3pt}=\hspace{-3pt} V_r^T\psi(x),\hspace{1.5pt} \bar{\psi}_{h}(x) \hspace{-3pt}=\hspace{-3pt} V_{n_L-r}^T\psi(x).
\end{align}
$W_{h\psi}$ has a full column rank $r>p$ and hence
\[
\psi_h(x_k) = (W_{h\psi}^TW_{h\psi})^{-1}W_{h\psi}^Th(x_k) \quad \forall \quad x_k \in \mathcal{M}.
\]
Therefore, under a similarity transformation $T \hspace{-4pt}= \begin{bmatrix}V_r & V_{n_L-r}\end{bmatrix}$, the model (\ref{eq: Koopman system with output}) takes the form of (\ref{eq: theorem 1 similarity transform}) with $T^{-1} = T^T$ where the first $r$ lifting functions of $\Tilde{\psi}$ given by $\psi_h$ lie in the span of $h$.
\hspace*{\fill}~\qed\par
\end{pf}

\begin{rem}
The $rank(W_h)\geq p$ only if $h(x)$ is linearly independent $\forall x \in \mathcal{M}$.
\end{rem}

We see that if $W_h$ is full column rank, the outputs can be determined completely by the OC-KO observables and the $W_\psi$ in (\ref{eq: output dynamic invertibility}) exists. In the event that $n_L<p$, the map $y=W_h\psi(x)$ is a map from a low dimensional space $\psi(x) \in \mathbb{R}^{n_L}$ to a high dimensional space $y \in \mathbb{R}^{p}$ and factor conjugacy is not defined for that case (in Section \ref{sec: Fusion using Koopman operators}). Hence, we project the outputs to a lower dimensional space $z= W_h^Ty \in \mathbb{R}^{n_L}$ to find a conjugate map. The following two corollaries illustrate how output functions can be used to identify all or a subset of Koopman observables (in Section \ref{sec: Fusion using Koopman operators}). 

\begin{cor} \label{corollary: outputs form Koopman operator}
If $rank(W_h)=n_L<p$, we can construct a diffeomorphic map between the states and projected outputs $z= W_h^Ty$. The $z$  dynamics  are given by
\[z_{k+1}= W_h^TW_hK(W_h^TW_h)^{-1}z_k.\] For $z = H(\psi(x)) = W_h^TW_h\psi(x), H:\mathbb{R}^{n_L}\rightarrow\mathbb{R}^{n_L}$, the dynamics of $\psi(x)$ and $z$ are diffeomorphic conjugate and the eigenfunctions that capture their dynamics have a bijective map using (\ref{eq: diffeomorphic conjugacy eigenfunction map}).
\hspace*{\fill}~\qed\par
\end{cor}


\begin{cor}\label{cor: rank(W_h) = n_L =p}
For $rank(W_h) = n_L =p$, $W_h$ becomes an invertible square matrix and the output dynamics in Corollary \ref{corollary: outputs form Koopman operator} simplifies to 
\[y_{k+1} = g(y_k)=W_hKW_h^{-1}y_k.\] 
Hence, the dynamics of $y$ and $\psi(x)$ are diffeomorphic conjugate.
\hspace*{\fill}~\qed\par
\end{cor}

The dynamics of the output can be constructed when $W_h$ is full column rank. The dynamics of the entire output or the projected output has a diffeomorphic conjugacy with the dynamics of the OC-KO observables $\psi(x)$ depending on whether $p=n_L$ or $p>n_L$ respectively under the map $y= W_h\psi(x)$. The solution to the dictionary of Koopman observables that capture the state dynamics are generally nonunique; they have infinite KOs for a given system (\ref{eq: NL_sys_without_output}). The fusion of the states and the outputs impose a constraint on the observables of the KOs. We explore this in the following corollary. 

\begin{cor}
Given a finite number of nontrivial output equations (\ref{eq: NL_output_equation}), let $\mathcal{K}_{f}$ and $\mathcal{K}_{f,h}$ denote the set of all KOs consistent with (\ref{eq: NL_sys_without_output}) and (\ref{eq: nonlinear system with output}) respectively. Any KO that satisfies (\ref{eq: nonlinear system with output}) also satisfies (\ref{eq: NL_sys_without_output}) $\Rightarrow \mathcal{K}_{f,h} \subset \mathcal{K}_{f}$. In the case of Corollary \ref{cor: rank(W_h) = n_L =p}, the output space captures the complete eigenfunction space of a KO that solves (\ref{eq: NL_sys_without_output}). But $\mathcal{K}_{f}$ contains more KOs whose eigenfunctions can be constructed by taking repeated product of the current eigenfunctions which cannot be spanned by the outputs $\Rightarrow \mathcal{K}_{f,h} \neq \mathcal{K}_{f}$. Hence $ \mathcal{K}_{f,h} \subset \mathcal{K}_{f}$.
\hspace*{\fill}~\qed\par
\end{cor}

We explicitly see that the output equations place a constraint on the OC-KO observables and that the dynamics of the outputs can be constructed when $rank(W_h)=n_L$. In the case \textit{(ii)} of Theorem \ref{theorem: outputs span all lifting functions} where $rank(W_h)<n_L$, using (\ref{eq: theorem 1 case 2}) to construct the output dynamics, we see that
\begin{equation}\label{eq: theorem1 leakage of psi_h}
    y_{k+1} = U_r\Sigma_rV_r^TKV_r\psi_h(x_{k}) +U_r\Sigma_rV_r^TKV_{n-r}\bar{\psi}_h(x_{k}) 
\end{equation}
where $\bar{\psi}_h(x_{k})$ is a leakage term that cannot be represented in the output space. Hence, the output cannot be constructed. To capture more information on $\bar{\psi}_h(x_{k})$, we need to construct the time-delay embedded outputs as the lifting functions \cite{brunton2017chaos,arbabi2017ergodic,kamb2020time}. This is a typical approach used to construct KORs when outputs partially measure the states \cite{brunton2017chaos,arbabi2017ergodic,kamb2020time, boskic2020koopman, fujii2019data, avila2020data}. In the case of pure output measurements, the time-delay embedded outputs capture the maximum dynamics that is observed by the outputs. But, when we fuse it with the state dynamics, there is no guarantee that that the output dynamics capture the entire state dynamics. Hence, we consider a more general case where the outputs capture a portion of the state dynamics. Since the OC-KO has the structure of a linear time invariant system, we invoke the observable decomposition theorem from \cite{hespanha2018linear} to separate out the dynamics of the observable lifted states that can be fused with the output dynamics.


\begin{thm}\label{theorem: observable constrained on output}
Given a nonlinear system (\ref{eq: nonlinear system with output}) with KO (\ref{eq: Koopman system with output}). Suppose $\psi(x)$ is the dictionary of Koopman observables for (\ref{eq: Koopman system with output}), not necessarily state-inclusive.

Then there exists a similarity transformation  $T$ and a projection matrix $W_\psi \in \mathbb{R}^{n_o \times (N+1)p}$ for some $N \in \mathbb{Z}_{\geq 0},n_o\leq n_L$ such that the dynamics of 
\begin{enumerate}[label=(\textit{\roman*})]
    \item a subset of the observables $\psi_o(x) \in \mathbb{R}^{n_o}$ of the transformed Koopman operator (under $T$)
    \begin{equation}\label{eq: Koopman system observable decomposition form}
    \begin{aligned}
        \begin{bmatrix} \psi_o(x_{k+1}) \\ \Bar{\psi}_{o}(x_{k+1})\end{bmatrix} &= \begin{bmatrix} K_o & 0 \\ K_{\Bar{o}} & K_{22}\end{bmatrix}
        \begin{bmatrix} \psi_o(x_{k}) \\ \Bar{\psi}_{o}(x_{k}) \end{bmatrix}\\ 
        y_k &=   \begin{bmatrix} W_{ho} & 0 \end{bmatrix} \begin{bmatrix} \psi_o(x) \\[-3pt] \Bar{\psi}_{o}(x)\end{bmatrix}
    \end{aligned}
    \end{equation}\vspace{-10pt}
\end{enumerate}
\begin{equation*}
        T^{-1}KT \hspace{-3.3pt}=\hspace{-3.3pt} \begin{bmatrix} K_o & 0 \\[-3pt] K_{\Bar{o}} & K_{22}\end{bmatrix}\hspace{-3pt}, W_hT \hspace{-3.3pt}=\hspace{-3.3pt} \begin{bmatrix} W_{ho} & 0\end{bmatrix} \hspace{-3pt}, T^{-1}\psi(x) \hspace{-3.3pt}=\hspace{-3.3pt} \begin{bmatrix} \psi_o(x) \\[-3pt] \Bar{\psi}_{o}(x)\end{bmatrix}
    \end{equation*}
\begin{enumerate}[resume,label=(\textit{\roman*})]
    \item and the projected time-delay embedded output
    \begin{equation}\label{eq: projected time-delay embedded output}
    z_k = W_\psi\begin{bmatrix} y_k & y_{k+1}& \cdots & y_{k+N} \end{bmatrix}^T
\end{equation}
\end{enumerate}
are diffeomorphically conjugate. 
\end{thm}
\begin{pf}
The Koopman operator representation (\ref{eq: Koopman system with output}) is a linear time invariant model. Using the observable decomposition theorem (Theorem 16.2 in \cite{hespanha2018linear}), there exists a similarity transformation $T$ that takes the system (\ref{eq: Koopman system with output}) to the form (\ref{eq: Koopman system observable decomposition form}) such that the subsystem
\begin{align*}
    \psi_o(x_{k+1}) &= K_o \psi_o(x_{k}) \triangleq g_1(\psi_o(x_{k}))\\
    y_k &= W_{ho}\psi_o(x_k) 
\end{align*}
is completely observable; $\psi_o(x) \in \mathbb{R}^{n_o}$ can be uniquely reconstructed from the outputs. Note that the system being observable is different from the observable functions in the context of Koopman. Then, there exists a $N\in \mathbb{Z}_{\geq 0}$ such that $Np \geq n_o$ and
\[
    \begin{bmatrix} y_k \\[-3pt] y_{k+1} \\[-3pt] \vdots \\[-3pt] y_{k+N}\end{bmatrix} \hspace{-2.5pt}= \hspace{-2.5pt}
    \begin{bmatrix} W_{ho}\\[-3pt] W_{ho}K_o \\[-3pt] \vdots \\[-3pt] W_{ho}K_o^N \end{bmatrix}\psi_o(x_{k}) \hspace{-2.5pt}=\hspace{-2.5pt}  \mathcal{O}\psi_o(x_k)
\]
where $\mathcal{O} \in \mathbb{R}^{Np \times n_o}$ has full column rank $n_o$. Then we can define a projection $W_\psi = \mathcal{O}^T$ to get (\ref{eq: projected time-delay embedded output}) such that 
\[z_k = \mathcal{O}^T\mathcal{O}\psi_o(x_k)\triangleq H(\psi_o(x))\] 
where $\mathcal{O}^T\mathcal{O} \in \mathbb{R}^{n_o \times n_o}$ is square invertible. The dynamics of $z$ is given by
\begin{align*}
    z_{k+1} &= \mathcal{O}^T\mathcal{O}\psi_o(x_{k+1}) = \mathcal{O}^T\mathcal{O}K_o\psi_o(x_k) \\
    &= \mathcal{O}^T\mathcal{O}K_o(\mathcal{O}^T\mathcal{O})^{-1}z_k \triangleq g_2(z_k)
\end{align*}
For the map $H:\mathbb{R}^{n_o} \rightarrow \mathbb{R}^{n_o}$, there exists an inverse map $H^{-1}:\mathbb{R}^{n_o} \rightarrow \mathbb{R}^{n_o}, H^{-1}(z) = (\mathcal{O}^T\mathcal{O})^{-1}z_{k}$ since $\mathcal{O}^T\mathcal{O}$ is full rank and invertible. Using the above, we can see that
\begin{align*}
    H\circ g_1(\psi_o(x_k)) &= g_2 \circ H(\psi_o(x_k)) = \mathcal{O}^T\mathcal{O}\psi_o(x_{k+1})\\
    g_1 \circ H^{-1}(z_k) &= H^{-1} \circ g_2(z_k) = (\mathcal{O}^T\mathcal{O})^{-1}z_{k+1}.
\end{align*}
Hence, the dynamics of $z$ and $\psi_o(x)$ are diffeomorphically conjugate. 
\hspace*{\fill}~\qed\par
\end{pf}

\begin{rem}
The original time-delay embedded output evolves in a high dimensional space $\mathbb{R}^{Np}$ when compared to $\psi(x_o) \in \mathbb{R}^{n_o}$ since we are seeking a sufficiently large set of output measurements $N$, such that $N \times p \geq n_o$. In that case, the factor conjugate map cannot be established similar to the case in Corollary \ref{corollary: outputs form Koopman operator}.  This observation motivates the construction of a time-delayed output embedding. 
\end{rem}

When outputs partially measure the states, we see that time-delay embedded outputs have a diffeomorphic map with a subspace of the lifting functions under a similarity transform and the dynamics evolving in the two spaces are diffeomorphically conjugate. 

\begin{cor}
When $n_o = n_L$, $z_k$ has a diffeomorphic map with the entire dictionary of observables $\psi(x_k)$. In this case, $z_k$ can constitute a Koopman observable basis such that it captures the dynamics of $\psi(x_k)$.
\end{cor}

\begin{cor}
The scenario where $N=0$ results in case \textit{(ii)} of Theorem \ref{theorem: outputs span all lifting functions} with $\psi_h(x_k) = \psi_o(x_k)$. This simplifies the output dynamics (\ref{eq: theorem1 leakage of psi_h}) to 
\begin{equation*}
    y_{k+1} = \begin{cases}
    W_{ho}K(W_{ho}^TW_{ho})^{-1}W_{ho}^Ty_k &, p>n_o\\
    W_{ho}KW_{ho}^{-1}y_k &, p=n_o
    \end{cases}
\end{equation*}

\end{cor}

We see that the OC-KO architecture can fuse the state and output dynamics even if the outputs do not capture the entire state dynamics. The above analysis shows that the OC-KO structure is such that the lifting functions capture the dynamics of the time-delay embedded outputs. This is an implicit constraint in the model structure of OC-KO. 


State-inclusive observables (\ref{eq: definition state inclusive observable function}) are useful since we can recover the dynamics by simply dropping the nonlinear observables. We show a sufficient condition for the existence of state-inclusive OC-KOs using a similar argument developed in \cite{yeung2018koopmanGramian}. We prove the following lemma which plays a crucial role in showing the invariance of the basis in the series expansions of analytic functions. 

\begin{lem} \label{lemma: invariant dictionary}
Given a dictionary of observable functions $\mathcal{D} = \{\psi_1(x),\psi_2(x),\hdots\}$ where $\psi_r(x)=\displaystyle\prod_{i=1}^n x_i^{p_{r,i}}$ with $p_{r,i} \in \mathbb{Z}_{\geq0}$, the product of functions from $\mathcal{D}$ lies in $\mathcal{D}$. 
\end{lem}
\begin{pf}
Consider two functions  $\psi_{\alpha}(x),\psi_{\beta}(x) \in \mathcal{D}$. Their product is given by
\begin{equation*}
    \begin{aligned} 
        \psi_{\alpha}(x)\psi_{\beta}(x) \hspace{-2pt}&=\hspace{-2pt} \displaystyle\prod_{i=1}^n x_i^{p_{\alpha,i}} \displaystyle\prod_{i=1}^n x_i^{p_{\beta,i}} = \displaystyle\prod_{i=1}^n x_i^{(p_{\alpha,i}+p_{\beta,i})}\in \mathcal{D}.
    \end{aligned}
\end{equation*}
Since the product of two functions lies in $\mathcal{D}$, the product of any number of functions from $\mathcal{D}$ also lies in $\mathcal{D}$ as can be seen by taking repeated product of functions.
\hspace*{\fill}~\qed\par
\end{pf}

If the monomial basis in the Taylor series expansion are propagated by one time step, we encounter a product of these monomials and Lemma \ref{lemma: invariant dictionary} defines a set that captures all these monomial functions and their products. We use the this lemma to build on the result from \cite{yeung2018koopmanGramian}, which shows the existence of state-inclusive KO for (\ref{eq: NL_sys_without_output}), to find state-inclusive OC-KOs for (\ref{eq: nonlinear system with output}).

\begin{prop} \label{theorem: existence of Kx and Wh}
Given the nonlinear system of the form (\ref{eq: nonlinear system with output}), if the functions $f$ and $h$ are real analytic on the open set $\mathcal{M}$, then there exists an OC-KO representation of the form (\ref{eq: Koopman system with output}) in the region of convergence of the Taylor series expansion of $f$ and $h$.
\end{prop}

\begin{pf}
Let us consider a dictionary of polynomial lifting functions 
$
\mathcal{D} = \{\psi_1(x),\psi_2(x),\hdots\}, \quad \psi_r(x)=\displaystyle\prod_{i=1}^n x_i^{p_{r,i}}
$
 where $p_{r,i} \in \mathbb{Z}^+$. Since $f$ is real analytic in $\mathcal{M}$, for any $x_0 \in \mathcal{M}$, there exists a Taylor series expansion centered about $x_0$ that converges to $f(x)$ for any neighborhood of $x_0$. Suppose $f(x) = \begin{bmatrix}f_1(x) & \cdots f_n(x)\end{bmatrix}^T$, the Taylor series expansion of $f$ about $x_0$ yields
\begin{equation*}
    f_i(x) \hspace{-2.5pt}=\hspace{-2.5pt} f_i(x_0) + \frac{\partial f_i}{\partial x}\Bigg{|}_{x_0}x + x^T \frac{\partial^2 f_i}{\partial x^2}\Bigg{|}_{x_0}x + \cdots \hspace{-2.5pt}=\hspace{-2.5pt} \sum_{j=1}^{\infty} c_{ij}\psi_j(x).
\end{equation*}
where each term lies in $\mathcal{D}$. Suppose $x_k= \begin{bmatrix}x_{k,1} & \cdots x_{k,n}\end{bmatrix}^T$  where $x_{k,i}$ indicates the $i^{th}$ state at discrete time index $k$. $x_k$ propagated by one time step yields a linear combination of functions in $\mathcal{D}$ as \[x_{k+1,i} = f_i(x) = \sum_{j=1}^{\infty} c_{ij}\psi_j(x_k).\]
To construct a linear system, we propagate each function on the right hand side $\psi_j(x_k)$ by one time step
\begin{equation*}
        \psi_j(x_{k+1}) = \displaystyle\prod_{i=1}^{n} x_{k+1,i}^{p_{j,i}} = \displaystyle\prod_{i=1}^{n} \Big{(}\sum_{j=1}^{\infty} c_{ij}\psi_j(x_k)\Big{)}^{p_{j,i}}
\end{equation*}
Since Lemma \ref{lemma: invariant dictionary} states that the product of any number of functions in $\mathcal{D}$ lies in $\mathcal{D}$, $\psi_j(x_{k+1}) = \sum_{r=1}^{\infty} k_{ir}\psi_r(x_k)$. Concatenating the expressions of $x_{k+1,j}$ and $\psi_j(x_{k+1})$ for all $j$, we get the Koopman representation $\psi(x_{k+1}) = K\psi(x_k)$.

Similarly, for the output equation, we can expand each function $h_i$ in $h(x) = \begin{bmatrix}h_1(x) & h_2(x) & \cdots h_m(x)\end{bmatrix}^T$ using the Taylor series expansion about $x=x_0$ to yield
\begin{equation*}
    \begin{aligned}
       h_i(x) \hspace{-2pt}&= h_i(x_0) \hspace{-2pt}+\hspace{-2pt} \frac{\partial h_i}{\partial x}\Bigg{|}_{x_0} \hspace{-9pt} x + \hspace{-2pt}x^T \frac{\partial^2 h_i}{\partial x^2}\Bigg{|}_{x_0}\hspace{-9pt}x + \cdots = \hspace{-4pt}\sum_{j=1}^{\infty} w_{ij}\psi_j(x)\\
       \Rightarrow y_k& = W_h\psi(x_k)
    \end{aligned}
\end{equation*}
Hence, if $f$ and $h$ of the nonlinear system (\ref{eq: nonlinear system with output}) are analytic in the open set $\mathcal{M}$, a state inclusive OC-KO of the form (\ref{eq: Koopman system with output}) exists for that system. 
\hspace*{\fill}~\qed\par
\end{pf}

We see that state-inclusive OC-KOs exist for nonlinear systems whose dynamics ($f$) and output ($h$) functions are real analytic. This is a sufficient condition and not a necessary one. In the next section, we explore how to identify OC-KOs from data using the DMD formulation.

\section{DMD with output constraints}\label{sec: oc-deepDMD}
The OC-KO identification involves fusing two datasets (states and outputs) to learn the dynamics of the nonlinear system (\ref{eq: nonlinear system with output}). The DMD algorithms typically identify KOs on one dataset. We introduce the more general OC-DMD formulation that incorporates the output constraint in the DMD problem to identify the OC-KOs. The OC-DMD formulation is
\begin{equation}\label{eq: general_optimization_problem}
    \begin{split}
        \min_{K,W_h,\psi} ||\psi(X_F)&-K\psi(X_P)||^2_F\\
    \text{such that \hspace{25pt}} Y_P &= W_h\psi(X_P) 
    \end{split}
\end{equation}
where $||.||_F$ is the Frobenius norm. 
The equality constraint is very stringent as the presence of output measurement noise could result in overfit models. Hence, we pose the weaker problem
\begin{align}\label{eq: direct oc-deepDMD}
\min_{\psi,K,W_h}\Bigg{|}\Bigg{|}
    \begin{bmatrix}
        \psi(X_F)\\
        Y_P
    \end{bmatrix} -
    \begin{bmatrix}
        K\\
        W_h
    \end{bmatrix}\psi(X_P)\Bigg{|}\Bigg{|}^2_F
\end{align}
which concurrently solves for $\psi(x),K$ and $W_h$. We refer to this as the \textit{direct OC-DMD} problem formulation. To explicitly show the effect of the output dynamics in learning OC-KOs, we also propose the \textit{sequential OC-DMD} problem formulation where the following sub-problems are solved sequentially:
\begin{subequations}\label{eq: sequential_optimization}
    \begin{align}
        \intertext{\textbf{1. Identification of Koopman dynamics:}}
        \min_{\psi_x,K_{1}}||\psi_x(X_F)&-K_{1}\psi_x(X_P)||^2_F \label{eq: sequential_optimization a}\\
        \intertext{\textbf{2. Output Parameterization:}}
        \min_{\varphi_y,W_{h1}}\Big{|}\Big{|}Y_P&-W_{h1}
        \begin{bmatrix}
        \psi_x(X_P)\\
        \varphi_y(X_P)
        \end{bmatrix}
        \Big{|}\Big{|}^2_F \label{eq: sequential_optimization b}\\
        \intertext{\textbf{3. Approximate Koopman Closure:}}
        \min_{\varphi_{xy},K_{2}}
        \Bigg{|}\Bigg{|}
        \begin{bmatrix}
        \varphi_y(X_F)\\
        \varphi_{xy}(X_F)
        \end{bmatrix} &- K_2
        \begin{bmatrix}
        \psi_x(X_P)\\
        \varphi_y(X_P)\\
        \varphi_{xy}(X_P)
        \end{bmatrix} \Bigg{|}\Bigg{|}^2_F.\label{eq: sequential_optimization c}
    \end{align}
\end{subequations}
This problem is called sequential OC-DMD, because it obtains a solution for the output-constrained Koopman learining problem as a sequence of optimization problems.  The solution generated by sequential OC-DMD, yields an OC-KO of the form:
\begin{align*}
    \psi(x_{k+1})& = K\psi(x_{k})\\ 
    y_{k} &= W_h\psi(x_{k})
    \end{align*}
    where 
    \begin{equation}\label{eq: sequential optimization model structure}
    \begin{aligned}
    \psi(x)& =
    \begin{bmatrix}
        x\\
        \varphi_x(x)\\
        \varphi_y(x)\\
        \varphi_{xy}(x)
    \end{bmatrix}, \,\, K= 
    \begin{bmatrix}
        K_{1} & K_{12} & 0 & 0 \\
        K_{21} & K_{22} & 0 & 0 \\
        K_{31} & K_{32} & K_{33} & K_{34} \\
        K_{41} & K_{42} & K_{43} & K_{44} 
    \end{bmatrix},\\
    K_1 &= \begin{bmatrix}
        K_{11} & K_{12} \\
        K_{21} & K_{22} 
    \end{bmatrix}, \,\,
    K_2= 
    \begin{bmatrix}
        K_{31} & K_{32} & K_{33} & K_{34} \\
        K_{41} & K_{42} & K_{43} & K_{44} 
    \end{bmatrix},\\
    W_{h} &= \begin{bmatrix}
        W_{h11} & W_{h12} & W_{h13} & 0
    \end{bmatrix} 
\end{aligned}
\end{equation}
where $\psi(x) = \begin{bmatrix}x^T &  \varphi^T(x)\end{bmatrix}^T$, $x \in  \mathcal{M} \subset \mathbb{R}^n$, $\varphi_x: \mathcal{M} \rightarrow \mathbb{R}^{n_x}$, $\varphi_y: \mathcal{M} \rightarrow \mathbb{R}^{n_y}$, $\varphi_{xy}: \mathcal{M} \rightarrow \mathbb{R}^{n_{xy}}$,   $n+n_x+n_y+n_{xy} = n_L$ and the output matrix  $W_{h1} = \begin{bmatrix}W_{h11} & W_{h12} & W_{h13} \end{bmatrix}$.

Sequential OC-DMD works to first solve for the KO of the state dynamics (\ref{eq: sequential_optimization a}), without accounting for any output measurements.  This represents the Koopman operator obtained from standard dynamic mode decomposition (or E-DMD) algorithms. The next step in sequential OC-DMD (\ref{eq: sequential_optimization b}) solves for the projection equation to parameterize output functions in terms of the existing Koopman observables $\psi_x(x)$, as well as any necessary additional output observables $\varphi_y(x)$ required to predict the output equation.  The last step (\ref{eq: sequential_optimization c}) then incorporates additional state-dependent observable dictionary functions $\varphi_{xy}(x)$ to guarantee closure of the new Koopman $\varphi_{y}(x)$ observables from step 2  (\ref{eq: sequential_optimization c}).

\subsection{Equivalence of Solution Spaces for Sequential and Direct OC-DMD} 

We now elucidate the relationship between the solution space of sequential and direct OC-DMD optimization problems. To do so, we make use of the following proposition regarding coordinate transforms which considers a manifold $M$ of dimension n ($dim(M)=n$). 



\begin{prop} [Proposition 2.18 from \cite{nijmeijer1990nonlinear}]\label{proposition coordinate chart from nijmeijer}
Suppose that $dim(M) = n$ and that $f_1,...,f_n,k\leq n$ are independent functions about $p\in M$. Then there exists a neighborhood $U$ about $p$, and functions $x_{k+1},...,x_n$ such that $(U,f_1,...,f_k,x_{k+1},...,x_n)$ is a coordinate chart.
\end{prop}

We use the coordinate transformation to separate out the observables that capture the state dynamics and the output equations minimally to show the mapping between the OC-KOs from both OC-DMD problems. We go ahead to prove the equivalency of the OC-DMD algorithms. 

\begin{thm} \label{theorem sequential deepDMD}
The optimization problems (\ref{eq: direct oc-deepDMD}) and (\ref{eq: sequential_optimization}) are equivalent in the minimal solution space where the observables are independent.
\end{thm}
\begin{pf}
To prove the equivalency of (\ref{eq: direct oc-deepDMD}) and (\ref{eq: sequential_optimization}), we need to show that for every solution of (\ref{eq: direct oc-deepDMD}) given by (\ref{eq: Koopman system with output}), there is a solution for (\ref{eq: sequential_optimization}) given by (\ref{eq: sequential optimization model structure}) and vice versa. It is easy to see that  (\ref{eq: sequential optimization model structure}) fits into the structure of (\ref{eq: Koopman system with output}) directly without any modification. Hence, we only need to prove the reverse. Suppose the OC-KO 
\begin{equation} \label{eq: Thm opt eqn 2}
    \begin{aligned}
        \begin{bmatrix}
            x_{k+1}\\
            \varphi_{y^*}(x_{k+1})
        \end{bmatrix}&= 
        \begin{bmatrix}
            \tilde{K}_{11} & \tilde{K}_{12}\\
            \tilde{K}_{21} & \tilde{K}_{22}
        \end{bmatrix}
        \begin{bmatrix}
            x_k\\
            \varphi_{y^*}(x_k)
        \end{bmatrix}\\
        y_{k} &=  \begin{bmatrix}
            \tilde{W}_{h1} & \tilde{W}_{h2}
        \end{bmatrix}
         \begin{bmatrix}
            x_{k}^T &        \varphi^T_{y^*}(x_{k})
        \end{bmatrix}^T
    \end{aligned}
\end{equation}
where $\varphi_{y^*}:\mathcal{M}\rightarrow \mathbb{R}^{n_{y^*}}$ is the minimal dictionary of observables that solves  (\ref{eq: direct oc-deepDMD}). Let us consider the minimal dictionary of observables to capture the state dynamics (\ref{eq: NL_sys_without_output}) be $\varphi_x(x) = T_{11} \varphi_{y^*}(x), \varphi_{x}:\mathcal{M}\rightarrow \mathbb{R}^{n_{x}}$ such that $n_x\leq n_{y^*}$. $n_x = n_{y^*}$ when $\varphi_{y^*}$ is the minimal dictionary that captures (\ref{eq: NL_sys_without_output}). Since the functions $\varphi_{y^*}(x)$ is the minimal OC-KO solution  $dim(\varphi_{y^*}(x)) = n_{y^*}$ with $n_x$  independent functions $\varphi_x(x) = \varphi_x(\varphi_{y^*}(x)) = T_{11} \varphi_{y^*}(x)$. Using Proposition  \ref{proposition coordinate chart from nijmeijer}, a coordinate transformation $\Tilde{T}_1 = \begin{bmatrix} \tilde{T}_{11} & \Tilde{T}_{12} \end{bmatrix}$ exists such that 
\[
\varphi_{y^*}(x) \hspace{-3pt}=\hspace{-3pt}  \begin{bmatrix} \tilde{T}_{11} & \Tilde{T}_{12} \end{bmatrix}\begin{bmatrix} \varphi_x(x)\\ \tilde{\varphi}_y(x) \end{bmatrix},\quad \begin{bmatrix} \varphi_x(x)\\ \tilde{\varphi}_y(x) \end{bmatrix} \hspace{-3pt}=\hspace{-3pt} \begin{bmatrix} T_{11} \\ T_{12} \end{bmatrix}\varphi_{y^*}(x)
\]
where $\Tilde{T}^{-1} = \begin{bmatrix} T_{11}^T & T_{12}^T \end{bmatrix}^T$. Since $\varphi_x(x)$ is sufficient to capture the KO for (\ref{eq: NL_sys_without_output}), $\tilde{K}_{12}\Tilde{T}_{12}=T_{11}\tilde{K}_{22}\Tilde{T}_{12} =0 $ and the transformed dynamics are 
\begin{align*}
    &\begin{bmatrix}
        x_{k+1}\\
        \varphi_{x}(x_{k+1})\\
        \Tilde{\varphi}_{y}(x_{k+1})
    \end{bmatrix}\hspace{-5.2pt}= \hspace{-5.2pt}
    \begin{bmatrix}
        \tilde{K}_{11} & \tilde{K}_{12}\Tilde{T}_{11} & 0\\
        T_{11}\tilde{K}_{21} & T_{11}\tilde{K}_{22}\Tilde{T}_{11} & 0\\
        T_{12}\tilde{K}_{21} & T_{12}\tilde{K}_{22}\Tilde{T}_{11} & T_{12}\tilde{K}_{22}\Tilde{T}_{12}
    \end{bmatrix}\hspace{-7pt}
    \begin{bmatrix}
        x_k\\
        \varphi_{x}(x_k)\\
        \Tilde{\varphi}_{y}(x_{k})
    \end{bmatrix}\\
     &y_{k} \hspace{-3pt}= \hspace{-4pt}  \begin{bmatrix}
        \tilde{W}_{h1} & \tilde{W}_{h2}\Tilde{T}_{11} & \tilde{W}_{h2}\Tilde{T}_{12}
    \end{bmatrix}
     \begin{bmatrix}
        x_{k}^T &        \varphi^T_{x}(x_{k}) & \Tilde{\varphi}_{y}(x_{k})
    \end{bmatrix}^T
\end{align*}
Suppose there exists functions $\varphi_y(x) = T_{21}\Tilde{\varphi}_{y}(x)$ $\forall$ $x\in \mathcal{M}$ where $\varphi_{y}:\mathcal{M}\rightarrow \mathbb{R}^{n_{y}}$, $n_y\leq (n_{y^*}-n_x)$ that in addition to the functions $x,\varphi_x$ are sufficient to capture the output equation, a similar procedure to the above can be adopted to result in the model structure (\ref{eq: sequential optimization model structure}).Therefore, the two optimization problems are equivalent. 
\hspace*{\fill}~\qed\par
\end{pf}

This theorem shows that if an OC-KO representation exactly captures state and output dynamics, without any redundancy in any of the dictionary functions, for every solution in sequential OC-DMD, we can find a solution in direct OC-DMD and vice versa \cite{boyd2004convex}.  We thus see that (\ref{eq: direct oc-deepDMD}) and (\ref{eq: sequential_optimization}) are equivalent optimization problems. The equivalency does not imply they are the same optimization problem because the objective function that they solve are different \cite{boyd2004convex}. 

Specifically, the model structure of a solution obtained from sequential OC-DMD  (\ref{eq: sequential optimization model structure}) is more sparse than the model structure of a solution obtained from direct OC-DMD (\ref{eq: Koopman system with output}). The model (\ref{eq: sequential optimization model structure}) explicitly reveals that the output dynamics constrain the OC-KO learning problem through the dictionary functions $\varphi_y$ and $\varphi_{xy}$. The advantage of direct OC-DMD problem is that it solves only one optimization problem as opposed to sequential OC-DMD which solves three. We shall compare the performances of these algorithms in Section \ref{sec: Simulation Results} with theoretical and numerical examples. 

\subsection{Coordinate Transformations of Standardization Routines on System State and Output Data}
A common practice in model identification problem is to scale the variables using standardization or normalization techniques to ensure uniform learning of all variables. Standardization of a scalar variable $\Tilde{x}$ yields
\begin{align}\label{eq: standardization}
    \Tilde{x}_{standardized} = \frac{\Tilde{x}- \mu(\Tilde{x})}{\sigma(\Tilde{x})}
\end{align}
where $\mu$ and $\sigma$ are the mean and standard deviation of $\Tilde{x}$. It is important to keep track of how such affine transformations modify the structure of OC-KO when comparing theoretical and practical results.

\begin{prop} \label{eq: theorem scaling}
Given the nonlinear system with output (\ref{eq: nonlinear system with output}) has a Koopman operator representation with observables given by
\begin{align}
    \label{eq: z1}
    \begin{bmatrix}
    x_{k+1}\\
    \varphi(x_{k+1})
    \end{bmatrix} &= 
    \begin{bmatrix}
    K_{11} & K_{12}\\
    K_{21} & K_{22}
    \end{bmatrix}
    \begin{bmatrix}
    x_{k}\\
    \varphi(x_{k})
    \end{bmatrix}\\
    y_k &= 
    \begin{bmatrix}
    W_{h1} & W_{h2}
    \end{bmatrix}
    \begin{bmatrix}
    x_{k}\\
    \varphi(x_{k})
    \end{bmatrix}
\end{align}
if the states and outputs undergo a bijective affine transformation $\Tilde{x} = Px+b$ and $\Tilde{y} = Qy+c$ where $P,Q$ are non-singular, then the state dynamics are transformed to 
\begin{align}\label{eq: transformed_state_dynamics}
    \Tilde{\psi}(\Tilde{x}_{k+1}) &= \Tilde{K}\Tilde{\psi}(\Tilde{x}_{k}) \nonumber\\
    \Tilde{\psi}(\Tilde{x}_{k}) &= 
    \begin{bmatrix}
    \Tilde{x}_{k+1} & \Tilde{\varphi}(\Tilde{x}_{k+1}) &    1
    \end{bmatrix}^T\\
    \Tilde{K} &= 
    \begin{bmatrix}
    PK_{11}P^{-1} & PK_{12} & (\mathbb{I}-PK_{11}P^{-1})b\\
    K_{21}P^{-1} & K_{22} & K_{21}P^{-1}b\\
    0 & 0 & 1
    \end{bmatrix} \nonumber
\end{align}
and the output equations become
\begin{align}\label{eq: transformed_output_equation}
    \Tilde{y_k} &= \Tilde{W}_h\Tilde{\psi}(\Tilde{x}_{k})\\
    \Tilde{W}_h &= 
    \begin{bmatrix}
    QW_{h1}P^{-1} & QW_{h2} & (QW_{h1}P^{-1}b+c)
    \end{bmatrix}.\nonumber
\end{align}
\end{prop}

\begin{pf}
When the state undergoes an affine transformation $\Tilde{x} = Px+b$, since the transformation is bijective($P^{-1}$ exists), the dynamics of the transformed state ($\Tilde{x}$) are given by substituting $x = P^{-1}\Tilde{x} - P^{-1}b$ in (\ref{eq: z1}) to yield
\begin{align*}
    \begin{bmatrix}
    P^{-1}\Tilde{x}_{k+1} \hspace{-2pt}-\hspace{-2pt} P^{-1}b\\
    \varphi(P^{-1}\Tilde{x}_{k+1} \hspace{-2pt}-\hspace{-2pt} P^{-1}b)
    \end{bmatrix} \hspace{-4pt}&=\hspace{-4pt} 
    \begin{bmatrix}
    K_{11} & K_{12}\\
    K_{21} & K_{22}
    \end{bmatrix}\hspace{-4pt}
    \begin{bmatrix}
    P^{-1}\Tilde{x}_{k} \hspace{-2pt}-\hspace{-2pt} P^{-1}b\\
    \varphi(P^{-1}\Tilde{x}_{k} \hspace{-2pt}-\hspace{-2pt} P^{-1}b)
    \end{bmatrix}.
\end{align*}
We define new observable functions $\Tilde{\varphi}(\Tilde{x}) \triangleq \varphi(P^{-1}\Tilde{x} - P^{-1}b)$ and the complete vector valued observable as $\Tilde{\psi}(\Tilde{x}) \triangleq \begin{bmatrix}\Tilde{x} & \Tilde{\varphi}(\Tilde{x}) & 1\end{bmatrix}^T$. By algebraic manipulation, we get the transformed state dynamic as given in (\ref{eq: transformed_state_dynamics}).

 
  If the output undergoes an affine transformation $\Tilde{y} = Qy+c$, we derive the transformed output in terms of the transformed state: 
\begin{align*}
    \Tilde{y}_k &= Q h(x_k) +c= Q
    \begin{bmatrix}
    W_{h1} & W_{h2}
    \end{bmatrix}
    \begin{bmatrix}
    P^{-1}\Tilde{x}_k - P^{-1}b\\
    \Tilde{\varphi}(\Tilde{x})
    \end{bmatrix} +c
\end{align*}
By simple algebraic manipulation, we end up with the affine transformed output equation(\ref{eq: transformed_output_equation}). 
\hspace*{\fill}~\qed\par
\end{pf}

We see that the bias in the affine transformation constrains an eigenvalue of the transformed OC-KO to be equal to 1. This is very important to track when using gradient descent based optimization algorithms to solve for OC-KOs because they identify approximate solutions and this could push the unit eigenvalue outside the unit circle making it unstable. To avoid numerical error in such algorithms, we should constrain the last row of the Koopman operator as in (\ref{eq: transformed_state_dynamics}). 

In both the OC-DMD problems (\ref{eq: direct oc-deepDMD}) and (\ref{eq: sequential_optimization}), the state-inclusive observables $\psi$ are considered as free variables. To learn the observables, we use the deepDMD formulation (\ref{eq: deep DMD formulation}) of representing $\psi$ as outputs of neural networks. When we incorporate the deepDMD formulation to solve the OC-DMD problems, we refer to them as OC-deepDMD algorithms and the identified OC-KOs as OC-deepDMD models. 

\section{Simulation Results}\label{sec: Simulation Results}
\begin{figure}[ht]
\centering
  \includegraphics[width=8.3cm,keepaspectratio]{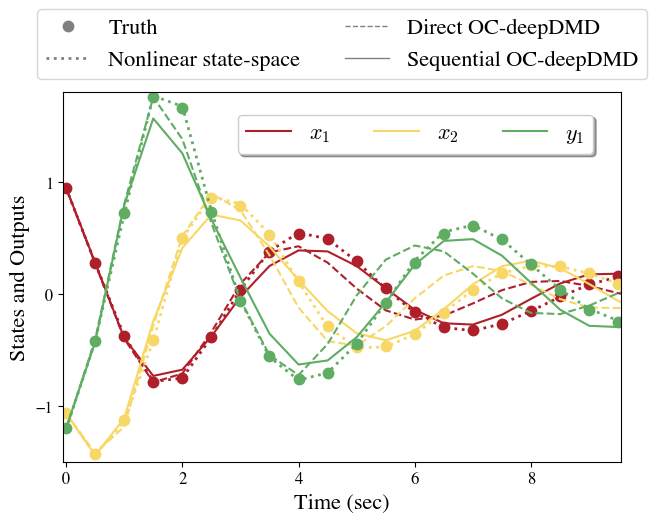}
  \caption{ \textit{MEMS Actuator:} Comparing the n-step predictions of states ($x_1,x_2$) and outputs ($y_1$) for the nonlinear state-space, direct OC-deepDMD and sequential OC-deepDMD models. }\label{fig: eg_MEMS_fit}
\end{figure}

\begin{figure}[ht]
\centering
  \includegraphics[width=8.4cm,keepaspectratio]{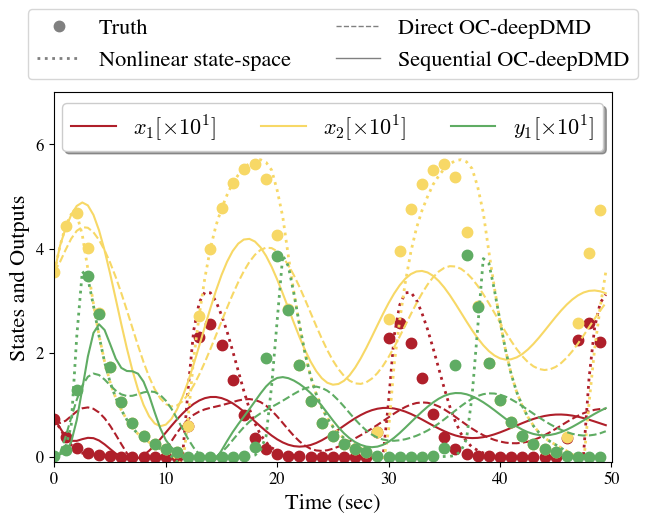}
  \caption{\textit{Activator-Repressor clock:} 
  n-step prediction comparisons of states ($x_1,x_2$) and outputs ($y_1$) for the nonlinear state-space and direct and sequential OC-deepDMD models. }\label{fig: eg_Act_Rep_fit}
\end{figure}

We consider three numerical examples in increasing order of complexity to evaluate the performance of the direct and sequential OC-deepDMD algorithms. The first example has an OC-KO with exact finite closure; there is a finite-dimensional basis in which the dynamics are linear. We use this as the benchmark for the comparison of the two algorithms. The other two examples do not possess finite exact closure. In those cases, we benchmark the proposed algorithms against nonlinear state-space models (with outputs) identified by solving 
\begin{align}\label{eq: Nonlinear state space model neural network optimization}
    \min_{f,h}||X_F - f(X_P)||_F^2 + ||Y_P - h(X_P)||_F^2
\end{align}
where the functions $f$ and $h$ are jointly represented by a single feed-forward neural network with $(n+p)$ outputs and we refer generally to this model, across multiple examples, as the \textit{nonlinear state-space model} (see captions in Figures \ref{fig: eg_MEMS_fit} and \ref{fig: eg_Act_Rep_fit}). 

\begin{table*}[th] 
\centering
\caption{Optimal hyper-parameters and performance of the oc-deepDMD models for all numerical examples}
\begin{tabular}{cccccc}
\toprule
System & Model & Hyperparameters & $r^2_{test}(x)$ & $r^2_{test}(x)$&  $r^2_{test}(y)$\\    
  & & & (1-step)  & (n-step) & \\
\midrule
    \multirow{4}{*}{ $\begin{array}{c}
         \text{Finite} \\
         \text{Closure} 
    \end{array}$}
 & \scriptsize{Deep DMD ($n_L=3$)} & \scriptsize{$n_x = 1,n_{xl} = 8,n_{xn} = 2$} & \scriptsize{$1$} & \scriptsize{$1$} & \scriptsize{-}
  \\
    & \scriptsize{Direct OC-deepDMD ($n_L=3$)} & \scriptsize{$n_x = 1,n_{xl} = 8,n_{xn} = 2$} & \scriptsize{$0.899$} & \scriptsize{$0.926$} & \scriptsize{$-0.147$}
  \\
    & \scriptsize{Direct OC-deepDMD ($n_L=5$)}& \scriptsize{$n_x = 3,n_{xl} = 7,n_{xn} = 5$} & \scriptsize{$1$} & \scriptsize{$1$} & \scriptsize{$1$}
    \\
    & $\begin{array}{c}
         \text{\scriptsize{Sequential OC-deepDMD}} \\
         \text{\scriptsize{($n_L=5$)}} 
    \end{array}$
     & \scriptsize{$\begin{array}{c}
        n_x = 1,n_{xl} = 8,n_{xn} = 2,n_y = 1, n_{yl} = 9, \\
        n_{yn} = 4, n_{xy} = 1, n_{xyl} = 7, n_{xyn} = 2
    \end{array}$} & \scriptsize{$1$} & \scriptsize{$1$} & \scriptsize{$1$}\\
  \midrule
    \multirow{3}{*}{ $\begin{array}{c}
         \text{MEMS} \\
         \text{Actuator} 
    \end{array}$}
 &\scriptsize{ Nonlinear state-space} & \scriptsize{$n_{xl} = 6,n_{xn} = 6$} & \scriptsize{$1$} & \scriptsize{$0.998$} & \scriptsize{$1$}
  \\
    & \scriptsize{Direct OC-deepDMD} & \scriptsize{$n_x = 6,n_{xl} = 3,n_{xn} = 6$} & \scriptsize{$1$} & \scriptsize{$0.84$} & \scriptsize{$0.995$}
    \\
    & \scriptsize{Sequential OC-deepDMD} & \scriptsize{$\begin{array}{c}
        n_x = 5,n_{xl} = 3,n_{xn} = 12,n_y = 1, n_{yl} = 8, \\
        n_{yn} = 6, n_{xy} = 3, n_{xyl} = 8, n_{xyn} = 3
    \end{array}$} & \scriptsize{$0.999$} & \scriptsize{$0.883$} & \scriptsize{$0.999$}\\
 \midrule
    \multirow{4}{*}{ $\begin{array}{c}
         \text{Activator-} \\
         \text{Repressor} \\
         \text{clock}
    \end{array}$}
 & \scriptsize{Nonlinear state-space} & \scriptsize{$n_{xl} = 5,n_{xn} = 6$} & \scriptsize{$1$} & \scriptsize{$0.854$} & \scriptsize{$0.999$}
  \\
    & \scriptsize{Direct OC-deepDMD} & \scriptsize{$n_x = 9,n_{xl} = 6,n_{xn} = 12$} & \scriptsize{$0.999$} & \scriptsize{$0.53$} & \scriptsize{$0.983$}
    \\
    & \scriptsize{Sequential OC-deepDMD} & \scriptsize{$\begin{array}{c}
        n_x = 3,n_{xl} = 9,n_{xn} = 8,n_y = 1, n_{yl} = 9, \\
        n_{yn} = 4, n_{xy} = 3, n_{xyl} = 9, n_{xyn} = 3
    \end{array}$} & \scriptsize{$1$} & \scriptsize{$0.349$} & \scriptsize{$0.998$}\\
    & $\begin{array}{c}
         \text{\scriptsize{Time-delay embedded }} \\
         \text{\scriptsize{Direct OC-deepDMD} }
    \end{array}$ & \scriptsize{$\begin{array}{c}
         n_x = 4,n_{xl} = 8, \\
         n_{xn} = 9, n_d = 6
    \end{array}$} & \scriptsize{$0.995$} & \scriptsize{$0.9116$} & \scriptsize{$0.8791$}\\
  \bottomrule
\end{tabular}\label{Table: All numerical examples}
\end{table*}

The neural networks in each optimization problem are constrained to have an equal number of nodes in each hidden layer. The hyperparameters for all the  optimization problems can be jointly given by $\{n_{ij}|i \in \{x,y,xy\},j \in \{o,l,n\}\}$ where $n_{io}$, $n_{il}$ and $n_{in}$ indicate the number of outputs, number of hidden layers and number of nodes in each hidden layer for the dictionary of observables indicated by $i$ ($\psi_x,\varphi_y$ or $\varphi_{xy}$). Sequential OC-deepDMD comprises $i \in \{x,y,xy\}$ $j \in \{o,l,n\}$, direct OC-deepDMD comprises $i \in \{x\}$   $j \in \{o,l,n\}$ and nonlinear state-space model comprises $i \in \{x\}$ and $j \in \{l,n\}$.

In each example, the simulated datasets are split equally between training, validation, and test data. For each algorithm, we learn models on the training data with various combinations of the hyperparameters. We train the models in Tensorflow using the Adagrad \cite{duchi2011adaptive} optimizer with exponential linear unit (ELU) activation functions. We use the validation data to identify the model with optimal hyperparameters for each optimization problem, which we report in 
Table \ref{Table: All numerical examples}. To quantify the performance of each model, we use the coefficient of determination ($r^2$) to evaluate the accuracy of the model predictions:
\begin{equation*}
    r^2 = 1 - \frac{||\Tilde{X} - \hat{\Tilde{X}}||_F^2}{||\Tilde{X}||_F^2} 
\end{equation*}
where $\Tilde{X}$ is the variable of interest ($X_F$ or $Y_P$) and $\hat{\Tilde{X}}$ is the prediction of that variable. We evaluate $r^2$ for the accuracy of 
\begin{itemize}
    \item \textit{the output prediction}:  $y_k= h(x_k)$ for the nonlinear state-space model and $y_k= W_h\psi(x_k)$ for the OC-deepDMD models
    \item \textit{1-step state prediction}: $x_{k+1} = f(x_k)$ for the nonlinear state-space model and $\psi(x_{k+1}) =K \psi(x_k)$ for the OC-deepDMD models. 
    \item \textit{n-step state prediction}:  $x_{i} = \underbrace{ f \circ f \circ \cdots \circ f}_{i \text{ times}}(x_0)$ for the nonlinear state-space model and $\psi(x_{i}) =K^i \psi(x_0)$ for the OC-deepDMD models where $x_0$ is the initial condition and $i$ is the prediction step.
\end{itemize}
The \textit{n-step state prediction} is a metric to test the invariance of the OC-KO; if the OC-KO is invariant, the $r^2$ for n-step predictions turns out to be 1. We do not consider the $n$-step output prediction as the error provided by that metric will be a combination of the errors in both state and output models.

\begin{exmp}\label{example theoretical closure}
    System with finite Koopman closure
\end{exmp}
Consider the following discrete time nonlinear system with an analytical finite-dimensional OC-KO \cite{brunton2016koopman}:
\begin{align*}
    \begin{bmatrix}
        x_{k+1,1}\\
        x_{k+1,2}
    \end{bmatrix} &=
    \begin{bmatrix}
        a_{11} & 0 \\
        a_{21} & a_{22} 
    \end{bmatrix}
    \begin{bmatrix}
        x_{k,1}\\
        x_{k,2}
    \end{bmatrix}+ 
    \begin{bmatrix}
        0\\
        \gamma x_{k,1}^2
    \end{bmatrix}
    \\
    y_k &= x_{k,1}x_{k,2}
\end{align*}
where $x_{k,i}$ and $y_k$ denote the $i^{th}$ state and the output at discrete time point $k$ respectively. We obtain the theoretical OC-KO using sequential OC-DMD (\ref{eq: sequential_optimization}):
\begin{itemize}
    \item \textit{Solving (\ref{eq: sequential_optimization a}) - } Adding the observable $\varphi_1(x) = x_1^2$ makes the dynamics linear:
    \begin{align}\label{eq: eg1 seq ocdeepDMD (a)}
        \begin{bmatrix}
            x_{k+1,1}\\
            x_{k+1,2}\\
            \varphi_1(x_{k+1})
        \end{bmatrix} &= 
        \begin{bmatrix}
            a_{11} & 0 & 0\\
            a_{21} & a_{22} & \gamma\\
            0 & 0 & a_{11}^2
        \end{bmatrix}
        \begin{bmatrix}
            x_{k,1}\\
            x_{k,2}\\
            \varphi_1(x_{k})
        \end{bmatrix}
    \end{align}
    \item \textit{Solving (\ref{eq: sequential_optimization b}) - } Adding $\varphi_2(x) = x_1x_2$ to  $\{x_1,x_2,\varphi_1(x)\}$ yields a linear output equation
    \begin{align*}
        y_{k} &= 
        \begin{bmatrix}
            0 & 0 & 0 & 1\\
        \end{bmatrix}
        \begin{bmatrix}
            x_{k,1}&
            x_{k,2}&
            \varphi_1(x_{k})&
            \varphi_2(x_{k})
        \end{bmatrix}^T.
    \end{align*}
    \item \textit{Solving (\ref{eq: sequential_optimization c}) - } To identify the dynamics of the added observable $\varphi_2(x)$ and ensure a closed basis, we add $\varphi_3(x) = x_1^3$ to get the OC-KO:
    \begin{align}\label{eg1 sim theoretical solution}
        \psi(x_{k+1})&= \begin{bmatrix}
            a_{11} & 0 & 0 & 0 & 0 \\
            a_{21} & a_{22} & \gamma & 0 & 0\\
            0 & 0 & a_{11}^2 & 0 & 0\\
            0 & 0 & a_{11}a_{21} & a_{11}a_{22} & a_{11}\gamma\\
            0 & 0 & 0 & 0 & a_{11}^3
        \end{bmatrix}\psi(x_{k}) \nonumber\\
        y_{k} &= \begin{bmatrix} 0 & 0 & 0 & 1 & 0\end{bmatrix}\psi(x_{k})
    \end{align}
    where $\psi(x_k)\hspace{-2pt} =\hspace{-2pt} \begin{bmatrix}
            x_{k,1} & x_{k,2} & \varphi_1(x_{k}) & \varphi_2(x_{k}) & \varphi_3(x_{k})
        \end{bmatrix}^T$.
\end{itemize}

\begin{figure}
  \includegraphics[width=8.4cm,keepaspectratio]{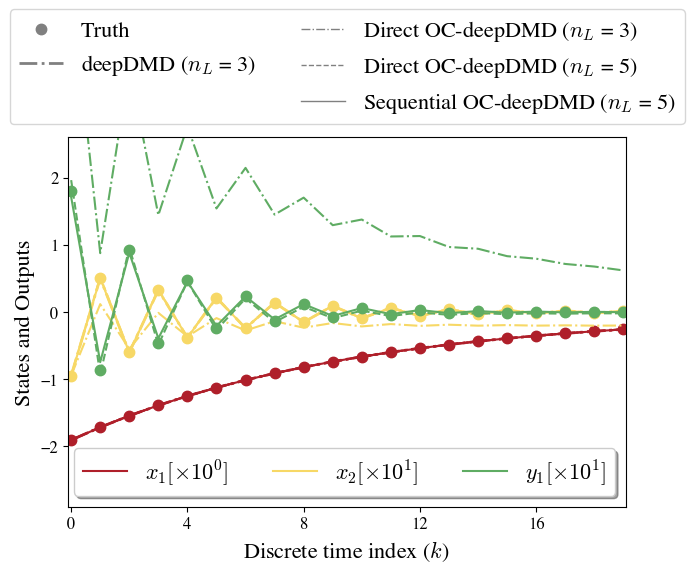}
  \caption{ \textit{Example \ref{example theoretical closure}:} Comparing the n-step predictions of states ($x_1,x_2$) and outputs ($y_1$) for the models: deepDMD with $n_L=3$, direct OC-deepDMD with $n_L=3,5$, and sequential OC-deepDMD with $n_L=5$ where $dim(\psi(x))=n_L$.
  }\label{fig: eg_TheoreticalExample_fit}
\end{figure}

We simulate the system to generate 300 trajectories, each with a different initial condition uniformly distributed in the range $5\leq x_{0,i} \leq 10,i=1,2$ and system parameters $a_{11} = 0.9, a_{21} = -0.4, a_{22} = -0.8, \gamma = -0.9$. The performance metrics of the identified models are given in Table \ref{Table: All numerical examples} and their $n-step$ predictions on a test set initial condition are shown in Fig. \ref{fig: eg_TheoreticalExample_fit}.

The deepDMD algorithm captures well the dynamics for $n_L=3$ and it matches with the theoretical solution (\ref{eq: eg1 seq ocdeepDMD (a)}). We learn the optimal direct OC-deepDMD model for $n_L=5$ Koopman dictionary observables. Fig. \ref{fig: eg_TheoreticalExample_fit} shows that the direct OC-deepDMD model with $n_L=3$ shows poor performance. Hence, we need additional observables to capture the output dynamics (given by $\varphi_2(x)$ and $\varphi_3(x)$ in (\ref{eg1 sim theoretical solution})), thereby validating Theorem  \ref{theorem: observable constrained on output}. We increase $n_L$ and identify both direct and sequential OC-deepDMD models. We observe that $n_L=5$ is the optimal value for both OC-deepDMD algorithms with $r^2\approx 1$ and agreeing with the theoretical solution (\ref{eg1 sim theoretical solution}).

\begin{figure*}[t]
  \includegraphics[width=\textwidth,keepaspectratio]{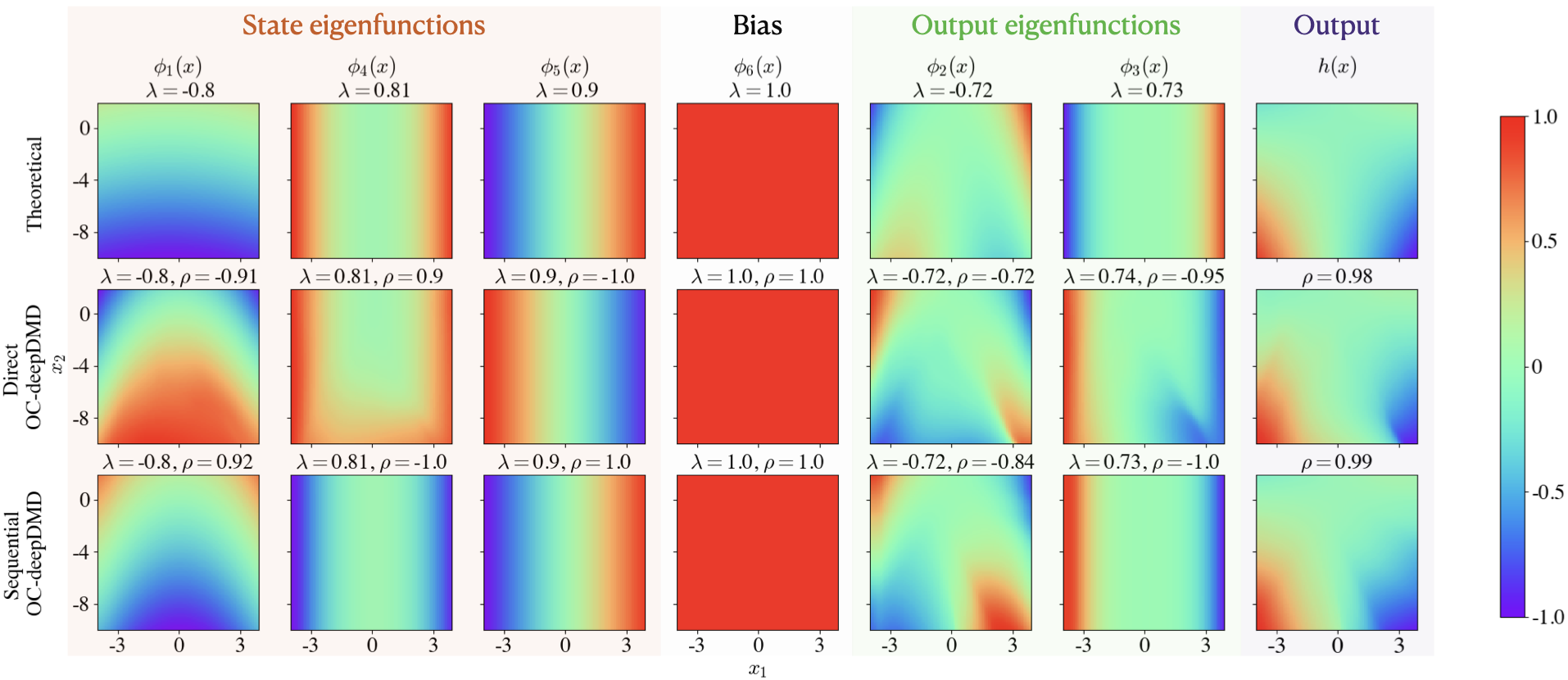}
  \caption{ \textit{Example \ref{example theoretical closure}:} The eigenvalues ($\lambda$) and corresponding eigenfunctions ($\phi$) computed from the theoretical (first row), direct OC-deepDMD (second row) and sequential OC-deepDMD (last row) models. The axes constitute the physical states of the system and the colorbar indicates the value of eigenfunctions normalized by the maximum absolute value attained by the corresponding eigenfunction. The Pearson correlation ($\rho$) is computed between $\phi$ of the OC-deepDMD and the theoretical models for each $\lambda$. Eigenfunctions $\phi_1, \phi_4, \phi_5$ with $\lambda=-0.8,0.81,0.9$ capture the state dynamics, additional $\phi_2, \phi_3$ with $\lambda = -0.72,0.73$ capture the output dynamics and $\lambda=1$ is due to the presence of the constant, unit-valued basis element (Proposition 3).}\label{fig: eg_TheoreticalExample_eig_func}
\end{figure*}

We evaluate the extent to which the OC-deepDMD algorithms capture the underlying system dynamics by comparing the eigenfunctions of the corresponding OC-deepDMD models with those of the theoretical OC-KO (\ref{eg1 sim theoretical solution}). Since the OC-deepDMD models are identified on standardized data, we use Theorem \ref{eq: theorem scaling} to reflect the transformation in the theoretical OC-KO. We compute the eigenfunctions of all the models using modal decomposition (\ref{eq: modal decomposition practical}). We observe that scaling and sign flip are two artifacts that lead to non-uniqueness of eigenfunctions as $\psi(x_{k+1})  = \sum_i v_i\lambda_i \phi_i(x)= \sum_i (-\alpha^{-1} v_i)\lambda_i (-\alpha\phi_i(x))$ where $\alpha$ is a nonzero scalar. We compensate for scaling by dividing each eigenfunction with its maximum absolute value to normalize it. When $r^2$ is computed under a sign flip, it leads to negative values. To account for the sign flip, we use the Pearson correlation ($\rho$) to compute the closeness between the normalized eigenfunctions.

We show the plot of the normalized eigenfunctions for the OC-deepDMD models and their correlation with the theoretical eigenfunctions in Fig. \ref{fig: eg_TheoreticalExample_eig_func}. We see that the sequential OC-deepDMD model captures both eigenvalues and eigenfunctions with a better accuracy than the direct OC-deepDMD model. This could be attributed to sequential OC-deepDMD model structure (\ref{eq: sequential optimization model structure}) being sparser among the two.  Sequential OC-deepDMD can explicitly track that the eigenfunctions corresponding to $\lambda = -0.8,0.81,0.9$ capture the state dynamics by solving (\ref{eq: sequential_optimization a}) and those corresponding to $\lambda = -0.72,0.73$ are added to capture the output dynamics by solving (\ref{eq: sequential_optimization b}) and (\ref{eq: sequential_optimization c}). This validates that the output dynamics do lie in a subspace of the OC-KO observables as proved in Theorem \ref{theorem: observable constrained on output}.

\begin{exmp} \label{example Nonlinear MEMS resonator with differential capacitor}
MEMS-actuator with a differential capacitor
\end{exmp}
\begin{figure}[ht]
\centering
  \includegraphics[width=5cm,keepaspectratio]{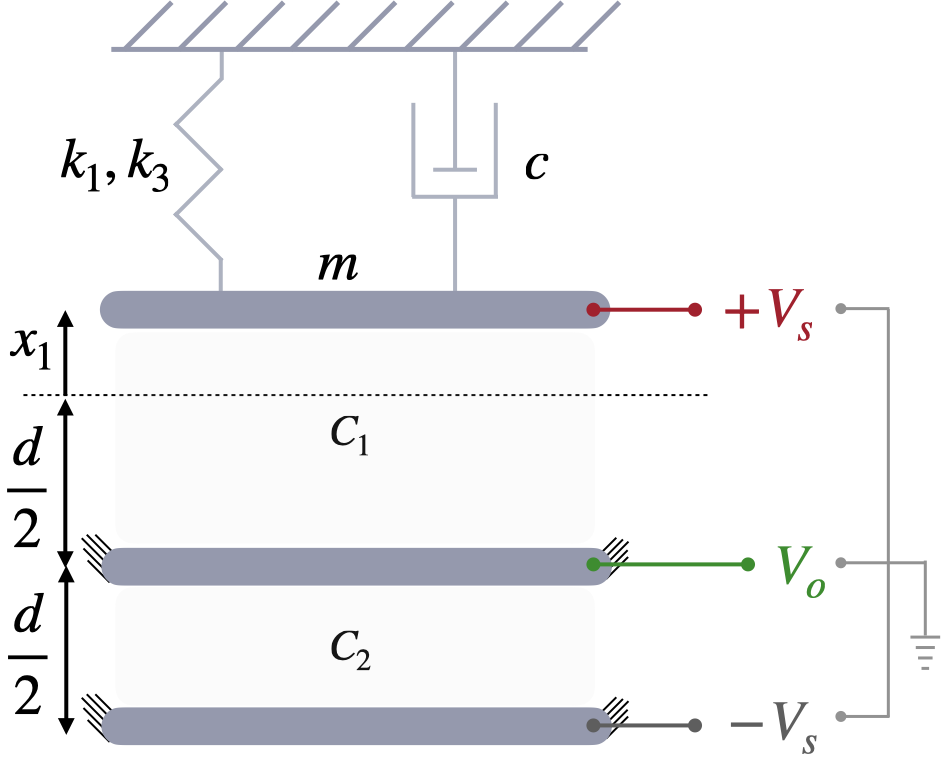}
  \caption{\textit{MEMS Actuator:} Schematic of a spring mass damper system as a MEMS actuator model. The displacement of the movable plate is measured by a differential capacitor with a fixed capacitance $C_2$ and variable capacitance $C_1$ by applying an input voltage $V_s$ and measuring the output voltage $V_o$. }\label{fig: eg_MEMSactuator}
\end{figure}
 We consider the free response of a MEMS resonator \cite{polunin2016characterization} modeled by a spring mass damper system with cubic nonlinear stiffness and a differential capacitive sensor to measure the displacement \cite{senturia2007microsystem} as shown in Fig. \ref{fig: eg_MEMSactuator}. It has the dynamics: 
\begin{align*}
    \dot{x_1} &= x_2\\
    \dot{x_2} &= -\frac{k_1}{m}x_1 -\frac{c}{m}x_2 -\frac{k_3}{m}x_1^3\\
    y &= V_o = \frac{C_1-C_2}{C_1+C_2}V_s = -\frac{x}{d+x}V_s
\end{align*}
where $x_1,x_2$ and $y$ are the displacement, velocity and output voltage measurements respectively. We simulate the system with the parameters $m=1,k_1 = 0.5,c = 1.0,k_3 = 1.0 \text{ and } V_{in}= 0.4$ to generate 300 trajectories with a simulation time of $15s$, a sampling time of $0.5s$ and initial condition, $x_0$, uniformly distributed in the range $(0,2)$. This system is more complex than Example \ref{example theoretical closure} as it has a single fixed point without a finite analytical OC-KO. Therefore, we benchmark the performance of the OC-deepDMD algorithms against the nonlinear state-space models identified by solving (\ref{eq: Nonlinear state space model neural network optimization}).

We see from Table \ref{Table: All numerical examples} that $r^2 \approx 1$ for 1-step state and output predictions of the nonlinear state-spacem, sequential, and direct OC-deepDMD models. When comparing the $n$-step predictions of these models on an initial condition from the test dataset (shown in Fig. \ref{fig: eg_MEMS_fit}), the nonlinear state-space model performs significantly better. Among the OC-KOs, the sequential OC-deepDMD performs marginally better (4\% higher accuracy). This indicates that all algorithms nearly accurately solve their respective objective functions which minimize the 1-step prediction error. However, the approximation of infinite-dimensional OC-KOs by finite observables lead to the OC-deepDMD models not being perfectly invariant. So, when the number of prediction steps increases, the error in the temporal evolution of the observables accumulates and propagates forward. A potential method to reduce the error in forecasting is to minimize multiple step prediction errors. We showcase this method prominently in the next example.

\begin{exmp} \label{example activator repressor 2 state}
Activator Repressor clock with a reporter
\end{exmp}
\begin{figure}[ht]
\centering
  \includegraphics[width=8cm,keepaspectratio]{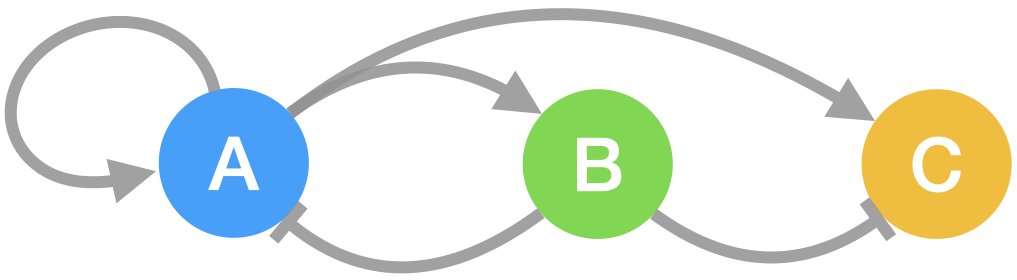}
  \caption{\textit{Activator-Repressor clock:} Schematic of the interaction of enzymes A, B and C (fluorescent reporter- output): (a) $A$ activates $A,B$ and $C$ (b) $B$ represses $A$ and $C$.}\label{fig: eg_activator_repressor}
\end{figure}

A more complex system is one with oscillatory dynamics that converge to an attractor. We consider the two state activator repressor clock \cite{del2015biomolecular} with the dynamics:
\begin{align*}
    \frac{dA}{dt} &= -\gamma_AA + \frac{\kappa_A}{\delta_A}  \frac{\alpha_A (A/K_A)^n + \alpha_{A0}}{1 +(A/K_A))^n + (B/K_B)^m} \\
    \frac{dB}{dt} &= -\gamma_B B + \frac{\kappa_B}{\delta_B} \frac{\alpha_B (A/K_A)^n + \alpha_{B0}}{1 +(A/K_A)^n}\\
    C &= \frac{(k_c/\gamma_c)A}{1 + (B/K_d)}
\end{align*}
where $A$, $B$ and $C$ are the conc. of enzymes with the network schematic as shown in Fig. \ref{fig: eg_activator_repressor}. $A$ and $B$ constitute the state and $C$ is the output fluorescent reporter assumed to be at steady state. We simulate the system using the parameters $\gamma_A=0.7$,  $\gamma_B = 0.5$, $\delta_A = 1.0$, $\delta_B = 1.0$, $\alpha_{A0}= 0.4$, $\alpha_{B0}= 0.004$, $\alpha_{A}= 0.2$, $\alpha_{B}= 0.2$, $K_A = 0.1$, $K_B=0.08$, $\kappa_A = 0.9$, $\kappa_B = 0.5$, $n=2$, $m=3$, $k_{3n}=3.0$  and  $k_{3d}=1.08$ to get an oscillatory behaviour. We generate 300 curves with a simulation time of $50s$, a sampling time of $0.5s$ and the initial conditions uniformly distributed in the interval $(0.1,1)$.

\begin{figure*}[t]
\centering
  \includegraphics[width=\textwidth,keepaspectratio]{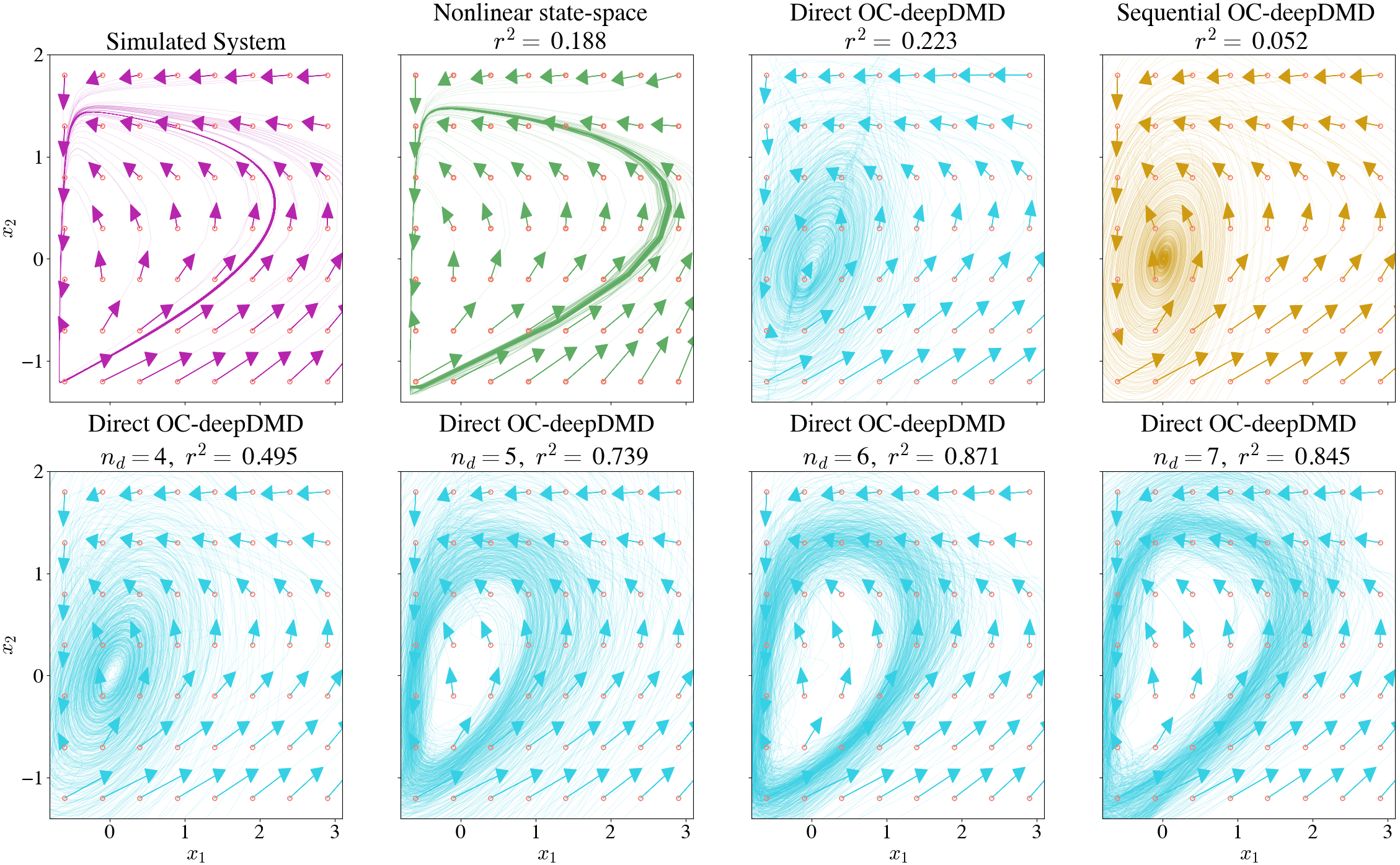}
  \caption{\textit{Activator-Repressor clock:} 
  The phase portraits of the theoretical, nonlinear state-space, direct OC-deepDMD and sequential OC-deepDMD models is shown in the first row. 
  The second row shows the phase portraits of the time-delay embedded OC-deepDMD models with the delay parameter $n_d=4,5,6$ and $7$. 
  The phase portraits are constructed using the same initial conditions highlighted by red dots and $r^2$ is computed using the theoretical phase portrait as the reference. }\label{fig: eg_activator_repressor_phase_plot}
\end{figure*}

We identify direct and sequential OC-deepDMD models and nonlinear state-space models and show the optimal model hyperparameters and $r^2$ values in Table \ref{Table: All numerical examples}. We see that $r^2 \approx 1$ for the 1-step and output predictions for all the models indicating that the corresponding algorithms nearly accurately solve their objective functions (which minimize 1-step prediction error) similar to the case in Example \ref{example Nonlinear MEMS resonator with differential capacitor}. The n-step predictions of the models in Fig. \ref{fig: eg_Act_Rep_fit} show that the nonlinear state-space model outperform both the OC-deepDMD models. But, here the direct OC-deepDMD model performs marginally better than the sequential  OC-deepDMD model (opposite of Example \ref{example Nonlinear MEMS resonator with differential capacitor}). Hence, we infer that both OC-deepDMD algorithms perform similarly. 

To evaluate how well the models capture the underlying dynamics, we construct phase portraits of the various models shown in Fig. \ref{fig: eg_activator_repressor_phase_plot}. We do so by considering initial conditions around the phase space of the limit cycle and plotting the n-step predictions of the models for each initial condition. We compare the phase portraits of each model with that of the simulated system using the $r^2$ metric. We see that the nonlinear state-space model captures a limit cycle but with an offset that results in a poor $r^2$ value. The OC-deepDMD models capture dissipating dynamics rather than that of a limit cycle. We speculate that the objective function is not sufficient to capture the dynamics and extend the objective function to minimize the error in multiple step predictions. To do so, we incorporate the idea we implement in \cite{balakrishnan2020prediction} to construct observables on the time-delay embedded states yielding the OC-KO  
\begin{align} \label{eq: causal jump oc-deepDMD model}
     \psi(x_{kn_d+n_d},\cdots,x_{kn_d+1}) &= K\psi(x_{kn_d},\cdots,x_{kn_d-n_d+1})\nonumber\\
     y_{kn_d}&= W_h\psi(x_{kn_d},\cdots,x_{kn_d-n_d+1}).
\end{align}
where $k$ is the discrete time index and $n_d$ indicates the number of time-delay embeddings. Since the two OC-deepDMD algorithms perform similarly, we stick to just using the direct OC-deepDMD algorithm to identify the time-delay embedded OC-deepDMD models. The phase portraits of the direct OC-deepDMD models as $n_d$ is increased is given in the second row of Fig. \ref{fig: eg_activator_repressor_phase_plot}.  We see that as $n_d$ increases, the phase portrait takes the structure of an oscillator with $n_d=6$ being optimal. 

We see from Table \ref{Table: All numerical examples} that the n-step prediction accuracy increases for this model at the expense of the 1-step and output predictions which reduce. This is because the formulation (\ref{eq: causal jump oc-deepDMD model}) simultaneously minimizes multiple step prediction errors \cite{balakrishnan2020prediction} which may not always yield optimal 1-step predictions. Hence, we see that OC-deepDMD algorithm has limitations when it comes to the case of oscillators and time-delay embedded OC-deepDMD models can be used to overcome them. 

\section{Conclusion}\label{sec: Conclusion}
In this work, we propose a novel method to fuse state and output measurements of nonlinear systems using Koopman operator representations that are augmented with a linear output equation (called OC-KOs). Using the concept of diffeomorphic conjugacy, we show that the dynamics of the measured output variables span a subspace of the OC-KO lifting functions and that the OC-KOs integrate the dynamics of both states and outputs. We show a sufficient condition for the existence of state-inclusive OC-KOs and propose two DMD algorithms that incorporate the output constraints to identify them. We use numerical examples to show the performance of these algorithms. 

In future work, we will use this technique to extract genotype-phenotype models of microbes by fusing their various time-series datasets. The genotype-phenotype models will enable us to control the persistence of these microbes in new environments. We expect the OC-KO based sensor fusion method to cater to a large range of dynamical systems where fusion of nonlinear dynamics of two measurement sets is desired for applications like observability analysis, observer synthesis and state estimation. 
\begin{ack}   
The authors would like to graciously thank Igor Mezic, Nathan Kutz, Jamiree Harrison, Joshua Elmore, Adam Deutschbauer, and Bassam Bamieh for insightful discussions.   
Any opinions, findings, conclusions, or recommendations expressed in this material are those of the authors and do not necessarily reflect the views of the Defense Advanced Research Project Agency, the Department of Defense, or the United States government.   

This work was also funded in part by the Department of Energy's Biological and Environmental Research office, under the DOE Scientific Focus Area: Secure Biosystems Design project, via gracious funding from Pacific Northwest National Laboratory subcontract numbers 545157 and 490521.
This work  was also partially supported by DARPA, AFRL under contract numbers FA8750-17-C-0229, HR001117C0092, HR001117C0094, DEAC0576RL01830, as well as funding from the Army Research Office's Young Investigator Program under grant number W911NF-20-1-0165.   Supplies for this work were partially supported by the Institute of Collaborative Biotechnologies, via grant W911NF-19-D-0001-0006.

\end{ack}


\bibliographystyle{unsrt}
\bibliography{main,sensor_fusion,sensor_list,system_identification,koopman_operator_sensor_fusion_applications,biology_sensor_fusion,Koopman_delay_embedding}
\end{document}